\documentclass[12pt]{article}
\normalbaselines
\font\aa=msam10 scaled\magstep1
\font\ab=msbm10 scaled\magstep1
\font\aba=msbm7 
\font\ac=eufm10 scaled\magstep1
\newcommand{\oo}{\openone}
\newcommand{\pree}{\aa \symbol{'004}}
\newcommand{\inee}{\ab \symbol{'044}}
\newcommand{\ine}{\mbox{\inee}}
\newcommand{\pre}{\mbox{\pree}}

\newcommand{\gz}{\mbox{\ab Z}}
\newcommand{\gr}{\mbox{\ab R}}
\newcommand{\gcp}{\mbox{\ab CP}}
\newcommand{\gc}{\mbox{\ab C}}
\newcommand{\gcm}{\mbox{\aba C}}

\newcommand{\Gras}{\mbox{$G_n({\gc}^{m+n})$}}
\newcommand{\got}[1]{{\mbox{\ac{#1}}}}
\newcommand{\mb}[1]{{\mbox{\boldmath{$#1$}}}}
\newtheorem{rem}{Remark}
\newtheorem{com}{Comment}
\newtheorem{lem}{Lemma}
\newtheorem{thm}{Theorem}

\newtheorem{pr}{Proposition}

\catcode `\@=11

\def\openone{\leavevmode\hbox{\small1\kern-3.8pt\normalsize1}}%

\def\@revmess#1#2{\typeout{REVTeX #1: #2}}

\newif\ifsecnumbers \global\secnumbersfalse
\def\ds@eqsecnum{\global\secnumberstrue}

\ifx\selectfont\undefined %
\@revmess{message}{NFSS not detected. Assuming OFSS.}
\fi
\let\reset@font\relax

\newcount\@indentflag \global\@indentflag=1 %
\newdimen\@eqtoeqnum \@eqtoeqnum=6pt %
\def\@indentamount{%
\ifcase\@indentflag 0pt\or\@centering\or0pt plus1fil\fi\relax
}
\def\FL{\global\@indentflag=0 }
\def\FR{\global\@indentflag=2 }

\newcounter{eqletter} 

\def\@eqnnum{\hbox{\reset@font\rm(\theequation)}}
\let\make@eqnnum=\@eqnnum %
\def\eqnum#1{\dec@eqnnum \global\def\make@eqnnum{\reset@font\rm(#1)}%

\def\@currentlabel{#1}%
}
\def\inc@eqnnum{\addtocounter{equation}{1}}
\def\dec@eqnnum{\addtocounter{equation}{-1}}
\@definecounter{equation}%
\@addtoreset{equation}{section} %
\def\theequation@prefix{\arabic{section}.} %

\def\mathletters{%
\inc@eqnnum  \setcounter{eqletter}{0}%
\edef\@currentlabel{\theequation}%
\def\theequation{\theequation@prefix\arabic{equation}.\alph{eqletter}}%
\def\inc@eqnnum{\addtocounter{eqletter}{1}}%
\def\dec@eqnnum{\addtocounter{eqletter}{-1}}%
}

\def\[{\relax
\ifmmode\@badmath\else\par\vskip-\lastskip\vskip\abovedisplayskip\fi
\hbox to\hsize\bgroup
\def\label##1{\@revmess{warning}{\string\label{##1} used
in \string\[\space environment after (\theequation)}}%
\hskip\@totalleftmargin\hskip\@indentamount$\displaystyle
}

\def\]{\relax
\ifmmode
$\hskip\@centering\egroup
\else
\@badmath
\fi
\vskip\belowdisplayskip
\global\@indentflag=1 %
\noindent\ignorespaces
}

\newbox\@testboxa
\newbox\@testboxb

\def\equation{\par\vskip-\lastskip\vskip\abovedisplayskip
\inc@eqnnum\let\@currentlabel=\theequation
\setbox\@testboxa=\hbox\bgroup\hskip\@totalleftmargin\hskip\@indentamount
\hbox\bgroup$\displaystyle
}

\def\endequation{$\egroup\hskip\@centering\egroup %
\setbox\@testboxb=\hbox{\make@eqnnum}%
\bgroup
\@tempdima\wd\@testboxa \advance\@tempdima by\wd\@testboxb
\ifcase\@indentflag
\advance\@tempdima by\@eqtoeqnum
\ifdim\@tempdima<\hsize %
\def\@tempa{0}%
\else
\def\@tempa{1}%
\fi
\or
\advance\@tempdima by2\@eqtoeqnum
\ifdim\@tempdima<\hsize %
\def\@tempa{0}%
\else %
\@tempdima\wd\@testboxa \advance\@tempdima by\wd\@testboxb
\advance\@tempdima by\@eqtoeqnum
\ifdim\@tempdima<\hsize %
\def\@tempa{0}%
\setbox\@testboxa\hbox{\hfill\box\@testboxa\kern\@eqtoeqnum}%
\else
\def\@tempa{1}%
\fi
\fi
\or
\advance\@tempdima by2\@eqtoeqnum
\ifdim\@tempdima<\hsize %
\def\@tempa{0}%
\setbox\@testboxb=\hbox{\kern\@eqtoeqnum\make@eqnnum}%
\else
\def\@tempa{1}%
3\fi
\fi
\ifnum\@tempa=0 %
\hbox to\hsize{\unhbox\@testboxa\box\@testboxb}%
\else %
\vbox{\hbox to\hsize{\unhbox\@testboxa}%
\vskip6pt %
\hbox to\hsize{\hfil\box\@testboxb}}%
\fi
\egroup
\global\let\make@eqnnum\@eqnnum %
\vskip\belowdisplayskip\noindent\global\@indentflag=1 \global\@ignoretrue
}

\def\eqnarray{\par\vskip-\lastskip\vskip\abovedisplayskip
\inc@eqnnum\let\@currentlabel=\theequation
\global\@eqnswtrue\m@th
\global\@eqcnt\z@
\tabskip\@totalleftmargin\advance\tabskip by\@indentamount\let\\\@eqncr
\halign to\hsize\bgroup\hskip\@centering
$\displaystyle\tabskip\z@{##{}}$&\global\@eqcnt\@ne
\hfil${{}##{}}$\hfil
&\global\@eqcnt\tw@ $\displaystyle\tabskip\z@{##}$\hfil
\tabskip\@centering \if@eqnsw\phantom{\make@eqnnum\kern\@eqtoeqnum}\fi
&\llap{##}\tabskip\z@\cr}

\def\endeqnarray{%
\@@eqncr\egroup
\vskip\belowdisplayskip\noindent
\dec@eqnnum\global\@indentflag=1
\global\let\make@eqnnum\@eqnnum %
\global\@ignoretrue
}

\def\nonumber{\global\@eqnswfalse
\def\label##1{\@revmess{error}{\string\label{##1} used
with \string\nonumber\space before (\theequation)}}%
}

\def\@@eqncr{\let\@tempa\relax
\ifcase\@eqcnt \def\@tempa{& & &}\or \def\@tempa{& &}%
\else \def\@tempa{&}\fi
\@tempa \if@eqnsw\make@eqnnum\inc@eqnnum\fi
\global\@eqnswtrue\global\@eqcnt\z@\global\let\make@eqnnum=\@eqnnum\cr
}

\@namedef{eqnarray*}{\def\@eqncr{\nonumber\@seqncr}%
\def\label##1{\@revmess{warning}{\string\label{##1} used
in eqnarray* environment before (\theequation)}}%
\eqnarray}
\@addtoreset{equation}{section}

\def\theequation{\arabic{section}.\arabic{equation}}

\def\section{\@startsection {section}{1}{\z@}{-3.5ex plus -1ex minus
     -.2ex}{2.3ex plus .2ex}{\normalsize\bf}}
\def\subsection{\@startsection{subsection}{2}{\z@}{-3.25ex plus -1ex minus
 -.2ex}{1.5ex plus .2ex}{\normalsize\bf}}

\def\rigla{
\vskip24pt
\hrule width\hsize\relax
\vskip 1.6cm
}

\def\@bibitem#1{\item\if@filesw \immediate\write\@auxout
{\string\bibcite{#1}{\the\value{\@listctr}}}\fi\ignorespaces}

\def\@cite#1#2{{#1\if@tempswa , #2\fi}} %
\def\@biblabel#1{$^{#1}$} %

\def\@lbibitem[#1]#2{\item[\@biblabel{#1}]\if@filesw
{\def\protect##1{\string ##1\space}\immediate
\write\@auxout{\string\bibcite{#2}{#1}}}\fi\ignorespaces}

\newdimen\WidestRefLabelThusFar
\global\WidestRefLabelThusFar\z@

\def\bibcite#1#2{\global\@namedef{b@#1}{#2}\@SetMaxRefLabel{#1}}

\def\@SetMaxRefLabel#1{%
{%
\setbox0\hbox{$^{\csname b@#1\endcsname}$}%
\ifdim\wd0>\WidestRefLabelThusFar
\global\WidestRefLabelThusFar=\wd0
\fi
}%
}

\def\@citex[#1]#2{%
\if@filesw\immediate\write\@auxout{\string\citation{#2}}\fi
\leavevmode\unskip$^{\scriptstyle\@cite{\@collapse{#2}}{#1}}$}

\def\CITE{%
\@ifnextchar[{\@tempswatrue\@CITEX}{\@tempswafalse\@CITEX[]}%
}

\let\onlinecite\CITE

\def\@CITEX[#1]#2{%
\if@filesw\immediate\write\@auxout{\string\citation{#2}}\fi
\leavevmode\unskip\ \@cite{\@collapse{#2}}{#1}}

\let\@bylinecite\cite

\def\@collapse#1{%
{%
\let\@temp\relax
\@tempcntb\@MM
\def\@citea{}%
\@for \@citeb:=#1\do{%
\@ifundefined{b@\@citeb}%
{\@temp\@citea{\bf ?}%
\@tempcntb\@MM\let\@temp\relax
\@warning{Citation `\@citeb ' on page \thepage\space undefined}%
}%
{\@tempcnta\@tempcntb \advance\@tempcnta\@ne
\edef\MyTemp{\csname b@\@citeb\endcsname}%
\def\@tempa{\@temptokena=\bgroup}%
\afterassignment\@tempa %
\@tempcntb=0\MyTemp\relax}%
\ifnum\@tempcntb=0\relax%
\@tempcntb=\@MM
\@citea\MyTemp
\let\@temp = \relax
\else %
\edef\@tempd{\number\@tempcntb}%
\ifnum\@tempcnta=\@tempcntb %
\ifx\@temp\relax %
\edef\@temp{\@citea\@tempd}%
\else
\edef\@temp{\hbox{--}\@tempd}%
\fi
\else %
\@temp\@citea\@tempd
\let\@temp\relax
\fi
\fi
}%
\def\@citea{, }%
}%
\@temp %
}%
}%

\catcode `\@=12

\begin{document}

\vspace*{2in}
\noindent
\begin{center}
{\bf ON THE GEOMETRY OF COMPLEX GRASSMANN MANIFOLD, ITS NONCOMPACT DUAL
 AND COHERENT STATES}
\vspace{1.3cm}\\
\end{center}
\noindent
\begin{center}
\begin{minipage}{13cm}
 S. Berceanu \\
Institute of Physics and Nuclear Engineering,
Department of Theoretical Physics, P. O. Box MG-6, Bucharest-Magurele, 
Romania;\\
 E-mail:  Berceanu@Roifa.IFA.Ro; Berceanu@Theor1.Ifa.Ro\\
\makebox[3mm]{ }\\
\end{minipage}
\end{center}

\vspace*{0.5cm}

\begin{abstract}
\noindent
Different topics on the differential geometry of the complex Grassmann manifold
are surveyed in relation to the coherent states. A calculation of the
 tangent
 conjugate locus and conjugate locus in the complex Grassmann manifold is
 presented. The proofs use the Jordan's stationary angles. Also
various formulas
for the distance on the complex Grassmann manifold  are furnished.

\end{abstract}

\vspace{1cm}

{\it AMS Mathematics Subject Classification: 53C22, 81R30\\

}

\newpage

\section{\hspace{-4mm}.\hspace{2mm}  INTRODUCTION}
\hspace*{0.8cm}
 
''Grassmann manifold... has been intensively studied for many years.  We
have not got a comprehensive knowledge of its geometry,
 however''.\cite{postnikov}

 Without entering into historical details, the
  Grassmann manifold has been intensively studied from the second half of
 the last century. The real  euclidean geometry of linear manifolds
 in a multidimensional
space was considered by Jordan\cite{jo} using only the methods of the
 analytic geometry.
  In the first half of our century the Grassmann
manifold was the main example in many constructions as the
 CW-cell decomposition,\cite{eh} the  Chern\cite{chern} and
 Pontrjagin\cite{pont} classes... The basic
 facts about the Grassmann manifold can be found in standard
 books.\cite{chern,bour}$^-$\cite{mil} Many recent references are
 based on the paper\cite{wong} of Y.-C. Wong. However, the modern reader 
has difficulties to follow Ref. \onlinecite{wong} because Wong uses some
notions as the stationary angles of Jordan\cite{jo} between two $n$-planes
 from an $n+m$
space. In fact, part of the results contained in the paper of Wong\cite{wong}
were known and they can be found in the papers of Rosenfel'd\cite{ros}
 and in his books.\cite{rosi1,rosi} 

In the case of the Grassmann manifold, the cut locus can be calculated 
explicitly.\cite{wong} The situation with the conjugate locus is more
 complicated.  Wong\cite{won} has published the expression of
 the conjugate locus in the Grassmann manifold and usually\cite{kn} his
paper is quoted as an example of a  calculation of
the conjugate locus in a  multidimensional manifold. The calculation of 
 Wong\cite{won} is essentially  based on a structure lemma the proof of which
 was
published later.\cite{wo} However, Wong has not published the proof of
his results on conjugate loci in Grassmann manifold. Sakai\cite{sak}
has calculated the tangent conjugate locus in the tangent space to the
 Grassmann manifold. He has observed that Wong's result announced in his
 paper\cite{won} is incomplete. Apparently,\cite{kob} this disagreement
of the results of Wong referring to the conjugate locus in the Grassmann
 manifold with the calculation of Sakai on the tangent conjugate locus has not
been  pointed out.

Among many other things, in this paper we present a proof of the results of
 Wong in
the complex Grassmann manifold and also another proof of the calculation of
Sakai in the tangent space to the Grassmann manifold. The part of the 
 conjugate locus calculated by Wong can be expressed as a Schubert variety.
The rest of the conjugate locus
  is characterized as  the subset of points of the
 Grassmann manifold which have at
least two of the stationary angles  equal. It contains as subset the subset of
 isoclinic spheres determined by Wong\cite{wwong} in connection with the
 Hurwitz\cite{hur} problem.

 The present paper can be considered from three points of view.

In this paper are put together many facts referring to the differential 
 geometry of the
 complex Grassmann
manifold. From this side, the paper has a survey character. However,   all the
proofs are original. The still open problem  refers to the conjugate locus,
as was already stressed. The proof uses also the stationary angles, which are
briefly presented. A short proof of the  structure  lemma of   Wong\cite{wo} is
given. Also explicit expressions for the distance on the
 Grassmann manifold  are deduced.

  On the other side, in this paper the geometry of the
complex Grassmann manifold is studied in relation to the coherent
 states.\cite{per,klauder}  The manifold of coherent vectors is the pull-back
of the dual of the tautological line bundle on a manifold identified with the 
$Det^*$ bundle\cite{chern} in the case of the complex Grassmann manifold. The
main observation is the fact that  the parameters which characterize the
 coherent
 states are in fact the Pontrjagin's coordinates of the $n$-plane. The proof of
the result of Wong on conjugate locus uses a parametrization that also appears
in the coherent state approach. 

This paper is a complete and self-contained  example of some notions related
 to
 the trial to find a geometrical characterization  of Perelomov's
construction of the coherent state manifold as K\"ahlerian embedding into a
projective space.\cite{sb}
 We remember that it has been pointed out that for
symmetric spaces the cut locus is equal to the polar divisor.\cite{sb1}
 This situation
is illustrated in the case of the complex Grassmann manifold. Also the
equality between the dimension of the projective space in which the Pl\"ucker
embedding takes place, the  Euler-Poincar\'e characteristic of the manifold and
the maximal number of orthogonal coherent vectors is true at least for flag
 manifolds.

The considerations below concern the geometry of the finite dimensional complex
Grassmann manifold \Gras, also denoted 

\begin{equation}
X_c=G_c/K=SU(n+m)/S(U(n)\times U(m))~. \label{xc}
\end{equation}

Simultaneously, some of the considerations
are made  also in the case of the noncompact dual of the compact Grassmann
manifold
\begin{equation}
X_n=G_n/K=SU(n,m)/S(U(n)\times U(m))~.\label{xn}
\end{equation}

  Most of the results are still true
for the infinite dimensional Grassmannian.\cite{ps}

The paper is organised as follows.

In \S 2 some basic facts about the complex Grassmann manifold are remembered.
 The Cauchy formula is 
still
true for projectively induced analytic line bundles over homogeneous K\"ahler
 manifolds. The Pontrjagin's coordinatization, the polar divisor  and cell
 structure are
considered in the Section 3. A rapid presentation of Schubert varieties
is proposed in Section 4 while the stationary angles are presented in \S 5.
The complex Grassmann manifold as symmetric space is treated in \S 6. In the 
same
Section is presented the connection between the Grassmann manifold and
the parametrization used in the coherent state approach.  The
 explicit expression of the exponential map which gives geodesics in the
 Grassmann
manifold is essential for calculating the conjugate locus in the manifold.
Lemma 6
will be used for determination of the tangent conjugate locus. The expression
of the diastasis function of Calabi,\cite{cal} recently used in the context of
coherent states,\cite{cgr}
 is given.  The cut locus
and conjugate locus are treated in \S 7. The main results are contained in
Theorem 2,  Proposition 2 and Comment 2. The last Section presents explicit
 expressions for the distance on the complex Grassmann manifold
(noncompact Grassmann manifold) which
 generalize the corresponding ones from the case of the Riemann sphere
(respectively, the disk $|z|<1$).

\hspace*{0.8cm}  

\section{\hspace{-4mm}.\hspace{2mm}   THE CAUCHY FORMULA}
\hspace*{0.8cm}2.1
Let us denote by $D_n({\bf K})$ the set of pure (decomposable) 
$n-$vectors of the exterior algebra $\bigwedge ^n {\bf K}$,
where {\bf K} is a complex vector space. For every $Z\in ~
D_n({\bf K})$, there exists  $n$ vectors 
 $\mb{z}_1,\ldots ,\mb{z}_n \in {\bf K}$ such that
\begin{equation}
Z=\mb{z}_1\wedge\ldots\wedge\mb{z}_n~. \label{z}
\end{equation}

The elements $ Z, Z'\in D_n({\bf K})$ are equivalent iff there exists
$\lambda \in {\gc}^{\star} \equiv {\gc}\setminus \{0\}$ such that
 $Z=\lambda Z'$,
that is, the associated $n-$subspaces $M_Z=~<\mb{z}_1,\ldots ,\mb{z}_n>$,
 $~M_{Z'}=\mbox {$<\mb{z'_1},\ldots ,\mb{z'_n}>$}$ are identical. So, there
is a canonical bijection of the set of $n-$subspaces of {\bf K}
on the image  $G_n({\bf K})$ of $D_n({\bf K})$ in the projective space
${\bf P}(\bigwedge ^n {\bf K})$ associated to $\bigwedge ^n {\bf K}$.
 $G_n({\bf K})$ is called the
Grassmannian of index $n$ of {\bf K}.
 The space of holomorphic sections of $
Det^*$ on  $G_n(K)$ is naturally isomorphic with $\bigwedge ^n ({\bf K}^*)$
(see Prop 2.9.2 in  Ref. \onlinecite{ps}). Here $E^* $ denotes the
 dual of the space $E$.

When ${\bf K}={\gc}^N,~ N= n+m$ and
 $\left\{\mb{e_i}\right\}_{i=1,\ldots ,N}$ is the canonical base
of {\bf K},  $G_n({\bf K})$ is the set of points of
 ${\bf P}(\bigwedge ^n {\bf K})$ with homogeneous (Pl\"ucker) coordinates
[$Z_I$]
\begin{equation}
Z=\sum _{I\subset{\cal S}(n,N)} Z_I{\mb e}_I ~.\label{2.2}
\end{equation}
If  $I_N=\{1,\ldots ,N\}$, then $\sigma:~I_N\rightarrow I_N$ is a
 Schubert symbol, i.e.
a permutation with the property that its restrictions to $I_n$ and
$I_N\setminus I_n$ are increasing and the set  of $N(n)={N\choose n}=
\displaystyle{\frac{N!}{n!m!}}$ 
Schubert symbols was denoted  by ${\cal S}(n,N)$.
 The Pl\"ucker embedding
$\iota :$ \Gras$\hookrightarrow {\gcp}^{N(n)-1} = {\bf PL},~{\bf L}= 
\bigwedge ^n{\gc}^N$,
\begin{equation}
\iota (Z)=[Z_I]_{I\subset{\cal S}(n,N)}\label{embpl}
\end{equation}
 is isometric and biholomorphic.\cite{leich} We have denoted
 $ [\omega ]\equiv \xi  (\omega )$,
 where  $\xi :{\bf K}\setminus \{0\}\rightarrow {\bf PK}$
 is the natural projection.

The $n-$vector $Z\not= 0$ is pure iff the (Pl\"ucker-)
Grassmann (-Cayley) relations are fulfilled, i.e.
\begin{equation}
\sum \epsilon_{i,J,H}Z_{J\setminus \{i\}}Z_{H\cup \{i\}}=0~,
\end{equation}
where $J, ~H\subset I_N,~ n=\#\{J\}-1=\#\{H\}+1$, and
 $\epsilon _{i,J,H}=+1~(-1)$ if the number of elements of $J$ and
$H$ less than $i$ have the same (resp. opposite) parity (cf.
Bourbaki\cite{bour}; see also   Ref. \onlinecite{ms} for the Hilbert
space Grassmannian).

Let $A^{i_1\ldots i_p}_{j_1\ldots j_p}$ denotes the minor of
order $p$ of the matrix $A$ whose elements are at the
intersection of the  rows $i_k$ with the columns $j_k,~
k=1,\ldots , p$. If in eq. (\ref{z}) 
\begin{equation}
\mb{z}_i=\sum_{a=1}^N \hat{Z}_{ia}\mb{e}_a~,\label{4}
\end{equation}
then
\begin{equation}
Z=\mb{z}_1\wedge\ldots\wedge\mb{z}_n=
\sum_{1\leq i_1<\ldots <i_n\leq N} Z^{i_1\ldots i_n}
\mb{e}_{i_1}\wedge\ldots\wedge\mb{e}_{i_n}	~,\label{zzz}
\end{equation}
where $ Z^{i_1\ldots i_n}=\hat{Z}^{1\ldots n}_{i_1\ldots i_n}$
are the Pl\"ucker coordinates  denoted in eq. (\ref{2.2}) by $Z_I$ and 
$\hat{Z}$ denotes the matrix $(\hat{Z}_{ia})_{1\le i\le n;1\le a\le N}$.

2.2 Let $((\cdot ,\cdot ))$ be the application $D_n({\bf K})\times
D_n({\bf K})\rightarrow \gc$ defined by the {\it Cauchy formula}
\begin{equation}
((Z',Z))\equiv (\iota (Z'),\iota (Z)) ~,\label{5.1}
\end{equation}
where $(\cdot ,\cdot )$ is the hermitian scalar product in ${\bf K}
\times {\bf K}$. The name of equation (\ref{5.1}) is justified by the {\it
Cauchy identity} (see eq. (6) p. 10 in Ref. \onlinecite{gant}) contained in

\begin{rem} \label{rem1}The following relation is true:
\begin{equation}
(\iota (Z'),\iota (Z))\equiv  
((\mb{z'}_1\wedge\ldots\wedge\mb{z'}_n,\mb{z}_1\wedge\ldots\wedge\mb{z}_n))=
\det [(\mb{z'_i},\mb{z}_j)]_{1\le i, j\le n}~.\label{6}
\end{equation}
\end{rem}
{\it Proof\/}:~ In fact
\begin{mathletters}\label{scz}
\begin{equation}\label{scaz}
\det [(\mb{z'_i},\mb{z}_j)]_{1\le i, j\le n}= 
\det(\hat{Z}'\hat{Z}^+), 
\end{equation}
\begin{equation}\label{scbz}
\det [(\mb{z'_i},\mb{z}_j)]_{1\le i, j\le n}= 
\det(\hat{Z}\hat{Z}'^+), 
\end{equation} 
\end{mathletters}
depending respectively on the convention of the hermitian scalar product
$(\cdot ,\cdot ):\newline~\mbox{{\bf K}$\times$ {\bf K}$\rightarrow$ {\gc}}$
\begin{mathletters}\label{sc}
\begin{equation}\label{sca}
(\mb{a},\lambda \mb{b})=\bar{\lambda} (\mb{a},\mb{b}),
\end{equation}
\begin{equation}\label{scb}
(\mb{a},\lambda \mb{b})=
\lambda (\mb{a},\mb{b}).
\end{equation}
\end{mathletters}
This corresponds respectively, to
\begin{mathletters}\label{9}
\begin{equation}
(\iota (Z'),\iota (Z))=\label{9a}\sum_{1\leq i_1<\ldots <i_n\leq N}
{Z}'^{i_1\ldots i_n}\bar{Z}^{i_1\ldots i_n}~,
\end{equation}
\begin{equation}
(\iota (Z'),\iota (Z))=\label{9b}\sum_{1\leq i_1<\ldots <i_n\leq N} 
\bar{Z}'^{i_1\ldots i_n} Z^{i_1\ldots i_n}~.
\end{equation}
\end{mathletters}
Eq. (\ref{6}) is a consequence of eqs. (\ref{scz}), (\ref{9}) and of  the
Binet-Cauchy formula: if
$A,~B, ~C$ are matrices  with $m\times n,~ n\times m$,
respectively $m\times m$ elements and $C=AB$, then (eq. (15) p. 9 in
 Ref. \onlinecite{gant})
\begin{equation}
\det C= \sum_{1\le k_1<\ldots <k_m\leq n}A^{1\ldots m}_{k_1\ldots k_m}
B^{k_1\ldots k_m}_{1\ldots m}~ .~\pre
\end{equation}

So, eqs. (\ref{5.1})-(\ref{scz}) justify the usual\cite{chern} definition of
 the {\it hermitian scalar product
of two pure $n-$vectors } ($n-${\it planes} of the Grassmannian), or,
 more precisely,
 of  the {\it hermitian scalar product in the holomorphic line bundle} $Det^*$:

\begin{equation}
((Z',Z))\equiv \det [(\mb{z'_i},\mb{z}_j)]_{1\le i, j\le n}~.\label{5.2}
\end{equation}
The infinite dimensional case can be found in Prop. 7.1 of
Ref. \onlinecite{ps}; see also eq. 2.10 in  Ref. \onlinecite{ms} .

2.3 If $Z, Z'\in G_n({\bf K})$, let $\theta$ be the angle
defined by the hermitian scalar product of two planes
\begin{equation}
\cos \theta (Z',Z)\equiv\frac{|((Z',Z))|}{\Vert Z'\Vert \Vert Z
\Vert }~.\label{10.1}
\end{equation}

Remark that $\theta$ in  equation  (\ref{10.1}) {\it is not the angle between
the two n-planes\/}, because $\theta $ is not invariant under
the motion group on the Grassmann manifold. The quantities
which are invariant under the group action are the $n$ stationary
angles $\theta_1,\ldots ,\theta_n$ of Jordan\cite{jo} related to 
$\theta$ by the relation  (\ref{55}). The only
situation in which the angle $\theta$ in relation (\ref{10.1}) is the angle
of the two $n-$planes occurs when the Grassmann manifold has
rank $1$, i.e. $r\equiv\min (m,n)=1$.

Eq. (\ref{5.1}) implies
\begin{equation}
\cos \theta (Z',Z)=\frac{|(\iota (Z'),\iota (Z))|}{\Vert \iota (Z')\Vert
 \Vert \iota (Z)
\Vert }~,\label{10.2}
\end{equation}
{\it and the r.h.s. of} eq. (\ref{10.2}) {\it defines}\cite{cay,coo}
{\it  the (intrinsic) distance on the geodesics joining} $\iota (Z'),
~\iota (Z)$ {\it in 
the
projective space} {\bf PL} {\it in which the Grassmann manifold is embedded},
\begin{equation}
\label{10.3.1} \cos d_c(\iota (Z'),\iota (Z))=
\frac{|(\iota (Z'),\iota (Z))|}{\Vert \iota (Z')\Vert \Vert \iota (Z)
\Vert }~.
\end{equation}

The elliptic hermitian distance, here called the Cayley distance,\cite{cay} is

\begin{equation}
d_c([\omega '],[\omega ])=
\arccos \frac{|(\omega ',\omega )|}{||\omega '||||\omega ||}~.
\label{cayley}
\end{equation}

The infinite dimensional case was treated by Kobayashi.\cite{kobi}

Now, it follows that

\begin{rem}[Rosenfel'd\cite{ros}] 
The angle  $\theta$\label{rem2}
 defined in eq. (\ref{10.1}) it is related to the Cayley distance
 $d_c$  by the relation
\begin{equation}
\theta (Z',Z)=d_c(\iota (Z'),\iota (Z))~.\label{10.4.1}
\end{equation}
\end{rem}

{\it Proof}: The Remark  results from eq. (\ref{10.2}) and eq.
 (\ref{10.3.1}).~\pre

Some authors (e.g. Study\cite{study}) prefer instead of the
 definition (\ref{cayley}) of
the distance $d_c$ the definition

\begin{equation}
d_c([\omega '],[\omega ])=
2\arccos \frac{|(\omega ',\omega )|}{||\omega '||||\omega ||}~,
\label{cayley1}
\end{equation}
which lead, instead of (\ref{10.4.1}) to
\begin{equation}
\theta (Z',Z)= 
\frac{1}{2}d_c(\iota (Z'),\iota (Z))\label{10.4.2}~.
\end{equation}
With the definition (\ref{cayley}), ((\ref{cayley1})) the elliptic  hermitian
 distance of two
points  on the Riemann sphere is one half the arc (respectively, the arc) of
 the great circle connecting the corresponding points of the Riemann
 sphere\cite{coo} (resp.\cite{study}) (see also \S 8). $d_c$ in eq.
 (\ref{cayley}) is equal to the minimum angle between the real lines
belonging to the complex lines (real 2-planes) in {\bf K} represented by
$[\omega ], [\omega ']\in \bf{PK}$.\cite{}

The Cauchy formula (\ref{5.1}) is still true\cite{sb,sb1} for  projectively 
induced\cite{hirz} analytic line bundles over homogeneous K\"ahler manifolds.

2.4 Now we briefly discuss the case of the noncompact  manifold $X_n$. 

Firstly, let us denote by
\begin{equation}
{\gcp}^{n-1,1}=SU(n,1)/S(U(n)\times U(1))\label{ncp}
\end{equation}
the hermitian hyperbolic space dual to ${\gcp}^{n-1}$. Then the noncompact
analogue of the distance (\ref{cayley}) is the hyperbolic hermitian Cayley
 distance

\begin{equation}
\label{ncay} d_c([\omega '],[\omega ])=
{\rm arccosh}\, \frac{|(\omega ',\omega )_n|}{||\omega '||_n||\omega ||_n}~,
\end{equation}
where the hermitian form on $\gc^N$, antilinear in the second entry
(convention (\ref{sca})), in the orthonormal basis, is

\begin{equation}
(\omega ',\omega )_n=\omega_1\bar{\omega '_1}-
\sum_{i=2}^n\omega_i\bar{\omega _i}~.\label{om}
\end{equation}

The noncompact manifold $X_n$  (\ref{xn}) admits the embedding
 $\iota ':X_n\hookrightarrow
{\gcp}^{N(n)-1,1}$. Eq. (\ref{5.1}) becomes
\begin{equation}
((Z',Z))_n\equiv (\iota '(Z'),\iota '(Z))_n ~,\label{55.1}
\end{equation}
and the {\it  Remark \ref{rem1}  with the r.h.s. in formula}
(\ref{6}) {\it replaced by}
 $\det [(\mb{z'_i},\mb{z}_j)_n]_{1\le i, j\le n}$ {\it is also
true  in the  case of the noncompact manifold $X_n$}.

The equation corresponding to eq. (\ref{10.1})  ((\ref{10.2})) is
\begin{equation}
{\rm \cosh}\, \theta (Z',Z)\equiv\frac{|((Z',Z))_n|}{\Vert Z'\Vert_n \Vert Z
\Vert _n}~,\label{10.1n}
\end{equation}
(respectively)
\begin{equation}
\cosh \theta (Z',Z)=\frac{|(\iota '(Z'),\iota '(Z))_n|}{\Vert 
\iota '(Z')\Vert_n \Vert \iota '(Z)
\Vert _n}~.\label{10.2n}
\end{equation}

So, {\it  Remark \ref{rem2} is also true in the noncompact case, with}
eq.  (\ref{10.1}) {\it replaced by}
(\ref{10.1n}) {\it and} eq. (\ref{cayley}) {\it replaced by} (\ref{ncay}).

\section{\hspace{-4mm}.\hspace{2mm} PONTRJAGIN'S  COORDINATIZATION, POLAR
DIVISOR AND CELLS}

\hspace*{0.8cm}  3.1 Let us consider $Z_0\in G_n({\bf K})$, where

\begin{equation}
Z_0=\mb{e}_1\wedge \ldots \wedge \mb{e}_n ~.\label{z0}
\end{equation}

Then we have the orthogonal decomposition 
\begin{equation}
{\bf K}= <Z_0>\oplus <Z^{\perp }_0>~,
\end{equation}
where $Z^{\perp }_0$ is the $m-$plane (completely) orthogonal to $Z_0$
defined by the $m$ vectors (\ref{24}).  Any $\mb{x}\in {\bf K}$ admits the
 decomposition

$$\mb{x}=\mb{u}\oplus \mb{v},~ \mb{u}\in <Z_0>, ~\mb{v}\in
<Z^{\perp }_0> ~,$$
and let $\varphi$ denotes the orthogonal projection $\mb{u}=\varphi
(\mb{x})$. More precisely, let us denote by $\varphi _{Z_0}(\mb{x})$ 
the orthogonal projection  $\mb{u}$ of the vector $\mb{x} \in {\bf K}$ on
the $n-$plane $Z_0$ in the direction $Z^{\perp }_0$.

Let us consider the (open) neighbourhood of the fixed $n-$plane $Z_0$
\begin{equation}
{\cal V}_{Z_0}=\{ Z\in G_n(\bf K) |~{\mbox {\rm projection~}} \varphi ~,
~\varphi _{Z_0}(Z)\subseteq Z_0 \mbox {{\rm~ is nondegenerate}}\} ~.\label{11}
\end{equation}

For a fixed $Z\in G_n({\bf K})$ let
\begin{equation}
\Sigma _Z=\{Y\in G_n({\bf K})~|~((Z,Y))=0\}~.\label{sigma}
\end{equation}
\begin{lem} \label{lem1}

  Let $Z_0,~ Z\in G_n({\bf K})$.  Then
$Z\in {\cal V}_{Z_0}$  iff one of the following equivalent conditions are
fulfilled

$$\begin{array}{llll}
 \mbox{ (A1) } &  ((Z,Z_0))\not=0, & &   \\
 \mbox{ (A2) } &  \varphi _Z(Z_0)= Z, & \mbox{or (A2')} &
 \varphi _{Z_0}(Z)=Z_0~,\\
 \mbox { (A3) } &  Z_0\cap \Sigma_Z=0, & \mbox{or (A3')} & Z\cap 
\Sigma_{Z_0}=0~.
\end{array}$$

 Equivalently, $Z \notin {\cal V}_{Z_0}$  iff one of the following
equivalent conditions are fulfilled
$$\begin{array}{ll}
  \mbox {(B1)}~~ &~~ ((Z,Z_0))=0~,\\
 \mbox{(B2)} ~~ & ~~\varphi _{Z}(Z_0)\subset Z,~ \varphi _Z(Z_0)\not= Z,
\mbox{ or (B2')}~~ \varphi_{Z_0}(Z)\subset Z_0, ~\varphi
_{Z_0}(Z)\not= Z_0 ~,  \\
 \mbox{(B3)}~~&~~Z_0\in \Sigma _{Z}~,~~~~~~~~~~~~~~~~~~~~\mbox{ or (B3')}~~Z\in
  \Sigma _{Z_0},
\end{array}$$
 and
\begin{equation}
\Sigma_{Z_0}=\{ Y\in G_n({\bf K})~|~\dim (Y\cap Z^\perp _0)\ge 1\}~.\label{13}
\end{equation}

 The complex Grassmann manifold can be represented as the  disjoint union
\begin{equation}
G_n({\bf K})={\cal V}_{Z_0}\cup \Sigma_{Z_0}~.
\end{equation}
\end{lem}
{\it Proof\/}: ~To prove  {\it (A)}, observe that the subspaces $Z,~Z'$
 {\it are related}
 (cf. Prop. 3.3 Ch. 7 in Ref. \onlinecite{hus}),
while $(B)$ can be obtained using the Remarks of Ch. I \S2 in  Ref.
\onlinecite{wu}, especially Lemma 1.3.
 See also Ch. 9 in Ref. \onlinecite{chern}. For the infinite dimensional 
Grassmannian see Ch. 7 of Ref. \onlinecite{ps}, especially Prop. 7.5.4.
~\pre

Geometrically, $\Sigma_{Z_0}$ is the cut locus of
$Z_0$, as was firstly observed by Wong\cite{wong} (also cf. Proposition 1
 below). The same property is true for a class
of spaces which generalizes the symmetric ones.\cite{sb,sb1}
$\Sigma_{Z_0}$ can be expressed as a Schubert variety  (cf. Lemma 3).

$\Sigma _X$ is called the {\it polar divisor} of $X$ (cf. Wu\cite{wu}).

Lemma  \ref{lem1} implies that for any $Z\in {\cal V}_{Z_0}$, there exists the
 vectors 
 $\mb{z}_1,\ldots ,\mb{z}_n \in \gc^N$ such that
relation (\ref{z}) holds, $\varphi _{Z_0}(Z)=
Z_0$ and $\varphi (\mb{z}_i)=\mb{e}_i$. Then, using the
 Pontrjagin\cite{pont} coordinates,

\begin{equation}
\mb{z}_i=\mb{e}_i+ \sum_{\alpha =n+1}^N Z_{i\alpha}\mb{e}_{\alpha},
~ i=1,\ldots ,n~,\label{15}
\end{equation}
and ${\cal V}_{Z_0}$ is  homeomorphic to $\gc^{n\times m}$.

 Let the vectors $\mb{z^\sigma_i}$
be such that $\hat{Z}^{\sigma} \in {\cal V}_\sigma $, where
\begin{equation}
\mb{z^\sigma _i}=\mb{e}_{\sigma (i)}+\sum_{\alpha =n+1}^N
Z_{\sigma (i)\sigma (\alpha )}\mb{e}_{\sigma (\alpha )}~,~i=1,\ldots ,n~ ,
\label{16}
\end{equation}
which for $\sigma$ identity was already given by eq. (\ref{15}). Then
 $(Z^{\sigma},{\cal V}_{\sigma }),~\sigma\in {\cal S}(n,N)$ furnish
an atlas of $G_n({\bf K})$, where $Z^\sigma=(Z_{\sigma (i),\sigma (\alpha
)})_{1\leq i\leq n,\, n+1\leq \alpha \leq n}$.

Let $\hat{Z}^{\sigma }$ be the (extended) matrix attached to the $n$ vectors
 (\ref{16}).
If $_\sigma Y$ denotes the submatrix of  $Y$
containing only the columns $\sigma (i), ~i=1,\ldots ,n$, then
\begin{equation}
_\sigma (\hat{Z}^\sigma )={\oo}_n,\label{nrx}
\end{equation} 
\begin{equation}
{\cal V}_{\sigma}=\{X\subset\Gras|\det\; {_\sigma}(\hat{X})\not=0\}.
\end{equation}

If $\hat{X}$ is the $n\times N$ matrix whose  $i-$th row consists of the
coordinates of the vectors ${\mb{x}}_i,~ i=1,\ldots ,n$, where $X={\mb{x}}_1
\wedge\cdots\wedge{\mb{x}}_n$, then
\begin{equation}
_\sigma \hat{X}= \hat{X}\Delta^\sigma ,\label{sigx}
\end{equation}
\begin{equation}
{\hat{X}}^\sigma =(_\sigma \hat{X})^{-1}\hat{X} ~,
\end{equation}
where
$$(\Delta^\sigma)_{ij}=\delta _{i\sigma (j)},~ i=1,\ldots ,N;~j=1,\ldots ,n.$$

The equations (\ref{nrx}) and (\ref{sigx}) imply that on ${\cal V}_\sigma
\cap{\cal V}_\tau\not=\emptyset$ a change of charts is given by the homographic
transformation  of the extended matrices
\begin{equation}
{\hat{Z}}^{\tau} =(\hat{Z}^{\sigma} \Delta ^\tau)^{-1}\hat{Z}^\sigma ,
~\sigma ,\tau
\in {\cal S}(n,N)\label{3xxx}.
\end{equation}
The  equations of the $n-$plane   $\hat{Z}^{\sigma }$ of 
$\gc^N$, generated by the $n$ vectors (\ref{16}),
  $\hat{Z}^{\sigma }\subset {\cal V}_{\sigma}$, are 
\begin{equation}
x_{\sigma(\alpha )}=\sum_{i=1}^n x_{\sigma (i)}
Z_{\sigma (i) \sigma (\alpha )}~,~ \alpha=n+1,\ldots ,N~,
\end{equation}
where $(x_1,\ldots ,x_N)$ are the local coordinates of ${\gc}^N$.

3.2 An ordering of the Schubert symbols is introduced as follows:
 $\sigma $ proceeds $\tau ~ (\sigma \prec \tau )$ if the least
 index $i,~ i\in I_n$ for which $\sigma (i)\not= \tau (i)$,
has the property  $\sigma (i)< \tau (i),$ where~ $ \sigma ,~\tau \in
{\cal S}(n,N)$.

Let
\begin{equation}
{\cal C}_{\sigma }=\{ \hat{Z}^{\sigma}\subset {\cal V}_{\sigma}~|~
{\det}( _{\tau }\hat{Z})=0,~ \sigma \prec \tau ,~ {\det} (_\sigma \hat{Z})
 \not= 0 \}~,
\end{equation}
and the matrix $\hat{Z}^\sigma$ is brought
to the  reduced echelon form:\cite{lu,mil} 

\begin{equation}
\hat{Z}^{\sigma }= \label{19}
\left(
\begin{array}{*{15}{l@{\hspace{-0.1pt}}}}
Z^1_1  & \ldots & Z^{\sigma (1)-1}_1 & r_1     & 0                  & 
\ldots & 0                  &    0 & 0  & \ldots  & 0 & 0 & 0 & \ldots 
 & 0 \\
Z^1_2  & \ldots & Z^{\sigma (1)-1}_2 & 0       & Z^{\sigma (1)+1}_2 &
\ldots & Z^{\sigma (2)-1}_2 & r_2 & 0  & \ldots  & 0 & 0  & 0  & \ldots
 & 0 \\
\vdots & \ddots & \vdots             & \vdots  & \vdots              &
\ddots & \vdots             & \vdots  & \vdots & \ddots  & \vdots & \vdots 
& \vdots  & \ddots & \vdots \\
Z^1_n  & \ldots & Z^{\sigma (1)-1}_n & 0       &  Z^{\sigma (1)+1}_n  &
\ldots & Z^{\sigma (2)-1}_n & 0 & Z^{\sigma (2)+1}_n & \ldots &
Z^{\sigma (n)-1}_n & r_n & 0 & \ldots & 0
\end{array}
\right)
\end{equation}
with the elements $r_i=1,~ i=1,\ldots ,n$.

 Let the notation
\begin{equation}
\omega (i)= \sigma (i)-i,~ i=1,\ldots ,n ~.\label{20}
\end{equation}

${\cal C}_{\sigma }$ is  homeomorphic to an (open) cell of complex dimension
\begin{equation}
d(\sigma )= \sum_{i=1}^n \omega (i) ~.\label{21}
\end{equation}

For the complex Grassmannian the groups of cell chains coincide, due to the
even dimension of the cells, with the groups of cycles, and, the group of
frontiers being trivial, the homology groups are isomorphic with the groups
of cell chains. 

Normalizing to one the row vectors in eq. (\ref{19}) such that the last
 element is
positive, the reduced echelon form\cite{mil,lu}
 are reobtained.
The open (closed) cells correspond to $r_i>0$, (respectively, $r_i\ge 0$).

In the Theorem below use is made of some notions referring to the coherent
states. The usual notation  will be remembered in \S 6.

\begin{thm}\label{th1}
 For the Grassmann manifold  \Gras \
 we have the equality of the following numbers:

1.  the maximal number of orthogonal coherent vectors;

2.  the number of critical points of the energy function $f_H$  
associated to a Hamiltonian $H$  which is a linear combination with
unequal coefficients of the generators of the Cartan algebra;

3.  the minimal dimension $N(n)$  appearing in the Kodaira
 (here Pl\"ucker)
 embedding $\iota : G_n({\gc}^N)\hookrightarrow {\gcp}^{N(n)-1}$;

4.   the Euler-Poincar\'e characteristic of the manifold,
$\chi ( G_n({\gc}^N))$;

5.   the number of Borel-Morse cells which appears in the CW-complex
decomposition of the Grassmannian;

6. the number of global sections in the holomorphic line bundle $Det^*$;

7. the dimension of the fundamental representation in the Borel-Weil-Bott
theorem.
\end{thm}
{\it Proof}: The theorem is proved with the theorems 1 and 2 in Ref.
\onlinecite{sbcag} particularized for the
Grassmann manifold and using the Cauchy formula.~\pre

 Theorem \ref{th1} is a particular case of a theorem true for  flag
 manifolds.\cite{sb}

3.3 If $Z=(Z_{i\alpha })_{1\leq i \leq n<\alpha \leq N}$ describes a $n-$plane
$Z\in{\cal V}_0$, then
 the extended matrix $\hat{Z}$  is 
\begin{equation}
\hat{Z}=({\oo}_n Z) ~,\label{exte}
\end{equation}
and the scalar product in eqs. (\ref{scz})
 can be written down, respectively, as
\begin{mathletters}
\begin{equation}
((Z',Z))=\det ({\oo}_n +Z'Z^+) ~,
\end{equation}
\begin{equation}
((Z',Z))=\det ({\oo}_n +ZZ'^+) ~.\label{23.2}
\end{equation}
\end{mathletters}
The noncompact analogue of eq. (\ref{23.2}) is
\begin{equation}
((Z',Z))_n=\det ({\oo}_n -ZZ'^+) ~.\label{23.2n}
\end{equation}

Given the $n-$plane $Z\in {\cal V}_0$ generated by the $n$ vectors in the
formula (\ref{15}), then the $m-$plane $Z^\perp $ orthogonal to $Z$ is
 generated by
\begin{equation}
\mb{z}^\perp_\alpha =\mb{e}_\alpha -\sum_{i=1}^n\bar{Z}_{i\alpha}
\mb{e}_i~,~\alpha=n+1,\ldots ,n+m~ ,\label{24}
\end{equation}
and
$$ (\mb{z}'^\perp_\alpha ,\mb{z}_i )=0~. \label{25}$$

Note also the following relations, corresponding
to the scalar product
(\ref{sca}) (respectively  (\ref{scb}))
\begin{mathletters}
\begin{equation}
((Z'^\perp ,Z^\perp ))=\det ({\oo}_m+Z'^+ Z)~,
\end{equation}
\begin{equation}
((Z'^\perp ,Z^\perp ))=\det ({\oo}_m+Z^+ Z')~.
\end{equation}
\end{mathletters}
Below we give a technical remark which has a clear geometrical meaning.

\begin{rem}\label{rem3} 
 If $Z,Z'\in {\cal V}_0\subset G_n({\bf K})$,
 then
\begin{equation}
((Z',Z))=\overline{((Z'^\perp ,Z^\perp ))}~,
\end{equation}
 or, explicitly,
\begin{equation}
\det ({\oo}_n+ZZ'^+)=\overline{\det}({\oo}_m+Z^+Z') ~.
\label{27.2}
\end{equation}

 Similarly, for $X_n$
\begin{equation}
\det ({\oo}_n-ZZ'^+)=\overline{\det}({\oo}_m-Z^+Z') ~.
\label{27.2n}
\end{equation}
\end{rem}
{\it Proof\/}:~ We present an algebraic proof of eq. (\ref{27.2}) for
 ${\bf K}={\gc}^N$.
This equation 
is a consequence  of the Schur formulas I, II (cf. Ref. \onlinecite{gant}
 p. 46). Let a
matrix  be partitioned in $4$ blocks, where the matrices $A$ and $D$ are 
non-singular. Then

\begin{equation}
\det \left(
\begin{array}{ll}
A & B \\
C & D
\end{array}
\right)
=\det (A) \det (D-CA^{-1}B)=\det (A-BD^{-1}C)\det (D) ~.\label{schur}
\end{equation}

 The Remark 3 follows taking the matrices $A,~ B,~ C,~ D$ as , respectively
${\oo}_n,~ Z,\\-Z'^+,~ {\oo}_m$.~\pre

The Remark 3 implies, via Remark 2, that
\begin{equation}
\cos \theta (Z',Z)=\cos \theta (Z'^\perp ,Z^\perp ) ~.\label{28}
\end{equation}
 Equation (\ref{28}) follows geometrically from Lemma 4 below and the fact 
that two 
$n-$planes and their orthogonal complements have the same invariants (cf.
 \S   48 p. 110 in Ref. \onlinecite{jo} ). 

\section{\hspace{-4mm}.\hspace{2mm}  SCHUBERT VARIETIES }

\hspace*{0.8cm} 
4.1 A Schubert variety $Z(\omega )$ associated with the monotone sequence

\begin{equation}
\omega=\{ 0\leq\omega (1)\leq \ldots \leq\omega (n)\leq m\},\label{29}
\end{equation}
is  the subset of the Grassmannian

\begin{equation}
Z(\omega )=\left\{ X\in G_n({\gc}^{n+m}) 
\vert \dim (X\cap {\gc}^{\sigma (i)})\geq i\right\}~ ,\label{30}
\end{equation}
where the sequences $\sigma$ and $\omega$ are related through the relation 
(\ref{20})
(see Pontrjagin\cite{pont}).

Usually (cf. Ch. XIV  p. 316 in Ref. \onlinecite{hodge} and Ref.
 \onlinecite{chern}), a nested sequence
of planes $V_i$ of dimension $\sigma (i),~ i=1,\ldots ,n$ is attached 
to the sequence  (\ref{29}), and these are the planes considered in the
 definition (\ref{30})
instead of ${\gc}^{\sigma (i)}$. However, because the definition of the
Schubert variety is independent of the concrete sequence of  planes $V_i$
modulo a congruence, it is enough to take $V_i={\gc}^{\sigma (i)}$
(cf. Pontrjagin\cite{pont}). This definition is also adopted by Milnor and 
Stasheff.\cite{mil}

Instead of considering the sequence (\ref{30}), it is enough to consider the 
sequence of
"jumps"\cite{pont,mil}
\begin{equation}
I_{\omega }=\left \{ 0=i_0<i_1< \ldots <i_{l-1}<i_{l}=n\right\}~,\label{31.1}
\end{equation}
where
\begin{equation}
\omega (i_h)<\omega (i_{h+1}),~ \omega (i)=\omega (i_{h-1}),~ i_{h-1}<i
\leq i_h,~ h=1,\ldots ,l. \label{31.2}
\end{equation}

Then
\begin{equation}
Z(\omega )=\left\{ X\in  G_n({\gc}^{n+m}) \vert \dim (X\cap {\gc}^
{\sigma (i_h)})\geq i_h,~  i_h\in I_{\omega}\right\} .\label{32}
\end{equation}

Let us consider the set of  elements in ``general position"\cite{pont}
 in $Z(\omega )$
(in fact, the subset of {\it generic elements}
in the sense of algebraic geometry\cite{hodge}):

\begin{equation}
Z'(\omega )=\left\{ X\in  G_n({\gc}^{n+m}) \vert \dim (X\cap {\gc}
^{\sigma (i_h)})= i_h,~  i_h\in I_{\omega}\right\} .\label{33}
\end{equation}

{\it The element  $Z\in G_n({\gc}^N)$,
$Z\in {\cal V}_0\cap
Z(\omega )\subset Z'(\omega )$, iff the coordinates in} eq. (\ref{15}) {\it 
verify
the condition}:\cite{pont,eh}

\begin{equation}
Z_{ij}=0,~~ j>\omega (i),~ i=1,\ldots ,n .\label{34}
\end{equation}

 Note that $Z'(\omega )$ defined by eq. (\ref{33}) corresponds to the open cell
 defined by eq. (8.15) in Chern's book,\cite{chern}
  and its complex dimension is given by equation (\ref{21}).

So, any $Z\in Z'(\omega )$ is locally characterized in ${\cal V}_0$
by a matrix $Z$ of the type

\begin{equation}
Z=\label{35}
\left( 
\begin{array}{ll}
Z_{i_1,\omega (i_1)} &{\bf 0}_{i_1,m-\omega (i_1)} \\
Z_{i_2-i_1,\omega (i_2)} & {\bf 0}_{i_2-i_1,m-\omega (i_2)} \\
\vdots & \vdots \\
Z_{i_l-i_{l-1},\omega (i_l)} & {\bf 0}_{i_l-i_{l-1},m-\omega (i_l)}\label{mat}
\end{array}
\right)
\end{equation}
and this representation makes very transparent the Chern\cite{chern}
 proof of the cell
decomposition  of the Grassmann manifold. In the formula above $Z_{p,q}$ 
denotes the $p\times q$ matrix and ${\bf 0}$ is the matrix with all elements
 $0$. The extended  matrix (\ref{19}) of the matrix (\ref{mat}) is obtained
using the  relations (\ref{16}).

4.2 Now the subset of the Grassmann manifold

\begin{equation}
V^p_l=\left \{Z\in  G_n({\gc }^{n+m})\vert \dim (Z\cap {\gc}^p)
\geq l\right\} ,\label{36}
\end{equation}
will be expressed as a Schubert variety, where $1<p< n+m$. Sometimes
another fixed $p-$plane of ${\gc}^N$, say $P_p$,  will be
considered in  eq. (\ref{36}) instead of ${\gc}^p$.
This situation will occur when  the Theorem 2 will be reformulated in the
notation of Wong,\cite{won} where $P_n=
{\bf O}$ and $P_m={\bf O}^\perp $.

Let also the notation
\begin{equation}
W^p_l=V^p_l-V^p_{l+1}=
\left \{Z\in  G_n({\gc}^{n+m})\vert \dim (Z\cap {\gc}^p)
= l\right\} ~.
\end{equation}

We shall prove the following structure Lemma of Wong\cite{won,wo} 

\begin{lem}\label{lema2} 
 Let $l\geq 1$  and $1<p<n+m$.  Then

\begin{equation}
V^p_l=Z(\omega ^p_l)~,\label{41}
\end{equation}
\begin{equation}
 W^p_l=Z'(\omega ^p_l)~.\label{42}
\end{equation}

 The following disjoint union is obtained

\begin{equation}
V^p_l=\label{43}
\left\{ 
\begin{array}{l}
\emptyset ~,{\mbox {\rm  if }} p<l~ {\mbox {\rm or }} l>n~
 {\mbox {\rm or }} p>l+m~, \\
G_n({\gc}^{n+m}),~ l=p-m~, \\
W^p_l\!\cup\! W^p_{l+1}\!\cup \!\cdots \!\cup\! W^p_{r_1-1}\!\cup\! W^p_{r_1},
~ \max (1,p\! -\! m\! +\!1)\! \leq \! l\!\leq \! r_1\! = \! \min (n,p).
 \end{array}
\right.
\end{equation}
 where
\begin{equation}
W^p_{r_1}=\label{44}
\left\{ \begin{array}{l@{\quad ~{\rm if}\quad}l}
{\gc}^n\subset G_n({\gc}^{n+m})~, & p=n~, \\
G_m({\gc}^{n+m-p})~, & p<n ~,\\
G_n({\gc}^p)~, & p>n ~.
\end{array}
\right.
\end{equation}
\end{lem}
{\it Proof\/}: ~The set of jumps (\ref{31.1}) for Schubert variety (\ref{36})
 is 
\begin{equation}
0=i_0<i_1=l<i_2=n ~. \label{jump}
\end{equation}

The relation $p=i_1+\omega (i_1)=\sigma (i_1)$ obtained forcing eqs.
 (\ref{36}) and (\ref{32}) to
coincide with the representations (\ref{34}), (\ref{35}) on the generic 
 elements (\ref{33}) imply that the sequence

\begin{equation}
\omega^p_l=(\underbrace{p-l,\ldots ,p-l}_l, \underbrace{m,\ldots ,m}_{n-l})
\end{equation}
is responsible for the Schubert variety (\ref{36}). The conditions that 
the variety (\ref{36})
to be nonvoid are $ m\geq p-l,~n-l\geq 0$. The case $p-l=m$ corresponds
to $\omega= (m,\ldots ,m)$ and then $V^{m+l}_l=G_n({\gc}^{n+m})$.

Eq. (\ref{jump}) implies that the matrix (\ref{35}) characterizing the set
(\ref{36}) has in this case only one submatrix
 ${\bf 0}_{l,m+l-p}$ with all elements zero. Then it  follows the disjoint
 union
\begin{equation}
V^p_l=W^p_l\cup V^p_{l+1}~.\label{39}
\end{equation}

In the representation (\ref{35}) the generic elements are characterized by
 the fact
that their first neighbours bordering the ${\bf 0}$ matrix in 
equation (\ref{35})
 are all nonzero, i.e.
\begin{equation}
\prod_{i=1}^l Z_{i\, p-l}\:\prod_{j=0}^{m+l-p}Z_{l+1\, p-l+j}\not=0 ~.
\end{equation}

The structure lemma is proved   iterating the splitting (\ref{39}) as far
 as possible.~\pre

Note that $V^p_l$ {\it is not} a (differentiable) (sub)manifold (of
 the Grassmannian).
In fact, $W^p_l$ consists of simple points of $V^p_l$ and $V^p_{l+1}$
is the singular locus (cf.  Ch. X \S 14
p. 87 in Ref. \onlinecite{hodge}) of $V^p_l$.\cite{won,wo} $V^p_l$ is an
 (irreducible) algebraic variety
of dimension  $l(p-l)+m(n-l)$, while $W^p_l$ is an (analytic) submanifold
of $G_n({\gc}^N)$ of the same dimension. 

In  particular, let $\Sigma _0$ be the polar divisor as defined by eq. 
(\ref{sigma})
of ${\bf O}\in G_n({\gc}^N)$. Then 

\begin{lem}[Wong,\cite{wong} Wu\cite{wu}] 
\label{lemma3}The polar divisor of the 
point {\bf O} is given  by
\begin{eqnarray}
\nonumber  \Sigma_0 & = & V^m_1 = Z(\omega^m_1)= Z(m-1,m,\ldots ,m)\\
& =  & \left\{ X\subset G_n({\gc}^{n+m})
\vert \dim (X\cap {\bf O}^{\perp})\geq 1\right\} .
\end{eqnarray}
\end{lem}

{\it Proof}:  The polar divisor is expressed as in eq. (\ref{13}).
 The relation (\ref{41}) implies the Lemma.~\pre
\section{\hspace{-4mm}.\hspace{2mm} THE STATIONARY ANGLES }
\hspace*{0.8cm}
5.1  Let $Z', Z$ be two $n-$planes of $G_n({\gc}^{n+m})$
given as in eq. (\ref{zzz}). Then the ($n$) {\it stationary angles} 
(see Jordan\cite{jo} for the real case),
{\it of which at most} $r=\min (m,n)$ {\it are nonzero}, are defined as the
{\it stationary} angles $\theta\in[0,\pi /2]$ between the vectors
\begin{equation}\label{abv}
\mb{a}=\sum_{i=1}^na_i\mb{z'_i},~\mb{b}=\sum_{i=1}^nb_i\mb{z}_i,~
\end{equation}
where
\begin{equation}
\cos \theta =\frac{|(\mb{a},\mb{b})|}{||\mb{a}||||\mb{b}||}~.\label{unghi}
\end{equation}
We shall prove a Lemma, which is implicitly contained in
 Jordan\cite{jo} :

\begin{lem}\label{lem4} 
 The squares
 $\cos ^2\theta _i$ {\it of the stationary angles between the\/ }
$\mbox{}n-${\it planes} $Z, Z'$ {\it with} $((Z,Z'))
\not= 0 $
 are given
as the eigenvalues of the matrix $W$
\begin{equation}
W=(\hat{Z}\hat{Z}^+)^{-1}(\hat{Z}\hat{Z'}^+)(\hat{Z'}\hat{Z'}^+)^{-1}
(\hat{Z'}\hat{Z}^+)~,\label{51}
\end{equation}
 which, for $Z, Z'\in {\cal V}_0$
 is
\begin{equation}
W=({\oo}+ZZ^+)^{-1}({\oo}+ZZ'^+)({\oo}+Z'Z'^+)^{-1}
({\oo}+Z'Z^+)~ .\label{53}
\end{equation}
\end{lem}
{\it Proof\/}: ~We adapt the Rosenfel'd's method (cf.
\S 3.3.15 p. 106  in Ref. \onlinecite{rosi1} and Ref. \onlinecite{rosi})
  to the complex Grassmann manifold in the Pontrjagin's coordinates.

Let us introduce the auxiliary function
\begin{equation}
U=(\mb{a},\mb{b})+\lambda (\mb{a},\mb{a}) +\mu (\mb{b},\mb{b}) ~.
\end{equation}

The calculation below is done with the condition (\ref{sca}) for the scalar
 product.

Taking the derivatives of $U$ with respect to $a_i$ and $b_i$, it follows that
\begin{equation}
\label{48}
\left\{ \begin{array}{l@{\;=\;0\; ,}}
(\mb{z'_i},\mb{z}_j)\bar{b}_j+
\lambda (\mb{z'_i},\mb{z'_j})\bar{a}_j\\
a_i(\mb{z'_i},\mb{z}_j)+
\mu b_i(\mb{z}_i,\mb{z}_j)
\end{array}
\right.
\end{equation}
which implies
\begin{equation}
\lambda \bar{\mu}=\cos^2\theta \geq 0.\label{49}
\end{equation}

Introducing the matrix of coordinates as in eq. (\ref{4}), eqs. (\ref{48}) can
 be written down in matricial form
\begin{equation}
\label{50}
\left\{ \begin{array}{l@{\;=\;0\; ,}}
\hat{Z'}\hat{Z}^+\bar{b}+\lambda\hat{Z'}\hat{Z'}^+\bar{a}\\
\hat{Z}\hat{Z'}^+\bar{a}+\bar{\mu}\hat{Z}\hat{Z}^+\bar{b}
\end{array}
\right.
\end{equation}
where $\hat{Z}, ~ \hat{Z'}$ are the extended $n\times N$ matrices (\ref{exte})
 attached to the $n-$planes $Z$, respectively $Z'$, and $a,~b$ are
$n-$column vectors.

It results from equations (\ref{50}) that $\bar{b}$ are the eigenvectors
corresponding to the eigenvalues $\lambda \bar{\mu}$ of the operator $W$ given
by eq. (\ref{51}) when $Z' \notin \Sigma _Z$, i.e. $((Z',Z))\not=0$.

Similarly, the scalar product with the convention $b)$ leads for the
 vector $b^t$ to the eigenvalues $\lambda \bar{\mu}$ of the operator
\begin{equation}
W^+=(\hat{Z}\hat{Z}^+)W(\hat{Z}\hat{Z}^+)^{-1}=
(\hat{Z}\hat{Z'}^+)(\hat{Z'}\hat{Z'}^+)^{-1}(\hat{Z'}\hat{Z}^+)
(\hat{Z}\hat{Z}^+)^{-1}~.
\end{equation}

The Lemma is proved taking into account eq. (\ref{49}).~\pre

Using the relations (\ref{48}), an algebraic proof of the theorems
 1-3 of Wong\cite{wong} follows.
A geometrical proof of these theorems in the case of the real Grassmann
manifold is given by Sommerville.\cite{som}

We shall show

\begin{lem}\label{lem5}
  Let $\theta$  be the angle defined by the 
hermitian scalar product in eqs. (\ref{10.1})-(\ref{10.4.1}),
 $d_c$  the Cayley distance
and  $\theta _1,\ldots ,\theta_n$  the stationary angles. Then
\begin{equation}
\cos \theta (Z,Z')= \cos d_c(\iota (Z'),\iota (Z))=
\cos \theta_1\cdots\cos \theta _n~.
\label{55}
\end{equation}
\end{lem}

{\it Proof}: It is  observed that Lemma \ref{lem4} implies
\begin{equation}
\det \sqrt{W}=\prod_{i=1}^n \cos\theta_i=
 {|\det ({\oo}+ZZ'^+) |
 \over |\det ({\oo}+ZZ^+)|^{1/2} |\det ({\oo}+Z'Z'^+)|^{1/2}} ~.
\label{54}
\end{equation}

But equations (\ref{5.2})-(\ref{10.4.1})  implies that 
$\cos \theta (Z,Z')
=\cos  d_c(\iota (Z'),\iota (Z))$ has also the expression (\ref{54}). ~\pre

Another proof of eq. (\ref{55}) can be found in Ref. \onlinecite{ros}
 or in more recent papers\cite{hi,fc}  which are based on the results of
 Wong.\cite{wong}

Now we attach an index $n$ to the $n$-plane $Z$ given by eq. (\ref{zzz}).

\begin{com} \label{com1}
 Let the $n' (n)$-plane $Z'_{n'}$ (resp. $Z_n$) with
 $n'\leq n$ such that $Z'_{n'}\cap Z_n=Z''_{n''}$. Then $n'-n''$ angles
of \/$Z'_{n'}$ and $Z_n$ are different from $0$ and $n''$ angles are $0$.
\end{com}

{\it Proof\/}:~ Firstly we consider the case $n'=n$. Then there are
 $n''$ common eigenvalues
of the matrix $W$ for which $\mb{a}$ and $\mb{b}$ are proportional. It
follows that  $\theta = 0$, corresponding to $\lambda\bar{\mu}= 1$
for $n''$ eigenvalues.

Let now $n'<n$. Then {\it there are at most $r_0=\min (n',n,m',m)$ angles which
 are different from} 0, where $n'+m'=n+m=N$. Observing that $\cos\theta =1$
iff $\mb{a}, ~\mb{b}$ are proportional, the considerations in the Comment are
still true even in this case.~ \pre

 The assertion contained in  Comment \ref{com1}  is largely discussed 
by  Sommerville in  Ch. IV p. 47 of Ref. \onlinecite{som} for
the general case for the real Grassmann manifold and also by Jordan
in Ref. \onlinecite{jo} \S  49 at p. 110. Reading the paper of Jordan, 
caution must be paid to the
fact that a $n$-plane in ${\gc}^{n+m}$ in Jordan's terminology is in fact
 a $m$-plane in the actual terminology.

5.2 We now briefly discuss the construction presented in this Section
 in the case of the noncompact manifold $X_n$.

If in eq. (\ref{unghi}) we consider instead of the hermitian scalar product
$(\cdot ,\cdot )$ the hermitian form
 $(\cdot ,\cdot )_n$ defined by (\ref{om}), then we could look for the
stationary ``angles'' (see also Wong\cite{wongg}) defined by the equation

\begin{equation}
{\rm cosh}\, \theta =\frac{|(\mb{a},\mb{b})_n|}{||\mb{a}||_n||\mb{b}||_n}~.
\label{unghin}
\end{equation}

With eq. (\ref{23.2n}), we find {\it the analogue of Lemma \ref{lem5}  in the
 case of the
noncompact manifold} $X_n$, $\cos ^2\theta_i$ {\it being substituted with} 
$\cosh ^2\theta_i$ {\it and} $W=W|_{\epsilon =-1}$, {\it where}
\begin{equation}
W(\epsilon )=({\oo}+\epsilon ZZ^+)^{-1}({\oo}+\epsilon ZZ'^+)
({\oo}+\epsilon Z'Z'^+)^{-1}({\oo}+\epsilon Z'Z^+)~ .\label{53n}
\end{equation}

 The noncompact analogue of eq. (\ref{54}) is
\begin{equation}
\det \sqrt{W(\epsilon )}=\prod_{i=1}^n {\rm cosh}\,\theta_i=
 {|\det ({\oo}+\epsilon ZZ'^+) |
 \over |\det ({\oo}+\epsilon ZZ^+)|^{1/2}
 |\det ({\oo}+\epsilon Z'Z'^+)|^{1/2}} ~.\label{54n}
\end{equation}
where $\epsilon  =-1$, while (\ref{55}) becomes

\begin{equation}
\cosh \theta (Z,Z')= \cosh d_c(\iota '(Z'),\iota '(Z'))=
\cosh \theta_1\cdots\cosh \theta _n~.\label{55n}
\end{equation}

\section{\hspace{-4mm}.\hspace{2mm}  THE COMPLEX GRASSMANNIAN AS
SYMMETRIC 
SPACE   AND COHERENT STATES}

\hspace*{0.8cm} 6.1. We  remember firstly the algebraic notation used in the
construction of symmetric spaces. 
 The Grassmann manifold is considered
 as compact hermitian irreducible Riemannian
globally  symmetric space of type A III.\cite{helg} We shall also remember the
relationship between the compact and noncompact Grassmann manifold. We use
 the conventions and notation from Ref. \onlinecite{sbl}.

$X_n$: the symmetric space of noncompact type (\ref{xn}).

$X_c$: compact dual of $X_n$ (\ref{xc}).

$o$: fixed base point of $X_n$ and $X_c$.

$K$: maximal compact subgroup of $G_n$, equal to the isotropy group of $G_n$
and $G_c$ at $o$.

$G^{\gcm}=G^{\gcm}_n=G^{\gcm}_c=G$: complexification of $G_c$ and $G_n$.

$\sigma$: Cartan involution, $\sigma=Ad\;s,~s=$ symmetry at $o$.

${\got g}_n,~{\got g},~{\got g}_c,{\got k}$: Lie algebras of $G_n,~ G,~G_c,~
K$, respectively.

${\got g}_n={\got k}+{\got m}_n$: sum of $+1$ and $-1$ eigenspaces of
 $d\sigma$.		

${\got g}={\got{g}}^{\gcm}_n={\got k}^{\gcm}+{\got m}_n$: complexification,
 where
$\got{m}={\got{m}^{\gcm}_n}$.
 
${\got g}_c={\got k}+{\got m}_c$: compact real form of \got{g},
 where $\got{m}_c=i{\got{m}_n}$.

{\got t}: Cartan algebra in ${\got g}_n$ and ${\got g}_c$.

${\got t}^{\gcm}$: Cartan subalgebra of {\got g}.

$\Delta$:  ${\got t}^{\gcm}$-root system of {\got g}, ${\got g}=
{\got t}^{\gcm}+
\Sigma_{\varphi\in\Delta}~{\got g}^{\varphi}$.

$\Delta_k$: set of compact roots, i.e. ${\got t}^{\gcm}$-root system of
 ${\got k}^{\gcm}$.

$\Delta_n$: set of noncompact roots, i.e. {\got m}-roots.

$\Delta^{\pm}$: set of positive and negative roots.

$(G_c,K)$: compact symmetric pair.

$({\got{g}}_c,d\sigma)$: orthogonal symmetric algebra of compact type
corresponding to $(G,K)$.

{\got a}: Cartan subalgebra of $(G_c,K)$, i.e. maximal commutative subset of
${\got m}_c,~ {\got a}\subset {\got t}$.

$\Sigma$: set of restricted roots, $\Sigma\equiv\Delta (\got{a};\got{g})=
 \{\gamma\in{\got a}|{\got g}^{\gcm}_{\gamma}\not=\{0\}\}$.

${\got g}^{\gcm}_{\gamma}$: restricted root space,
 i.e. ${\got g}^{\gcm}_{\gamma} =
\{X_{\gamma}\in{\got g}^{\gcm}|[H,X_{\gamma}]=\gamma
(H)X_{\gamma}, \forall H\in{\got a}\}$.

Let $\pi :G_c\rightarrow G_c/K$ the natural projection and $o=\pi (e)$, where
$e$ is the unity element in $G$. Then the subspace ${\got m}_c$ is naturally
identified  with the tangent space $(X_c)_0$ by means of the mapping $d\pi$.
The negative of the Killing form on ${\got g}_c$ defines an inner product $Q:
{\got m}_c\times {\got m}_c\rightarrow \gc$:
\begin{equation}
\label{qxy}
Q(X,Y)=-\frac{1}{2}{\rm{Tr}} (XY),
\end{equation}
which is $Ad\;G_c$ and $\sigma$-invariant. The Riemann
structure is defined by restricting $Q$ to ${\got m}_c\times {\got m}_c$ and
 translating with $G_c$. The geodesic map $\gamma_X:t\rightarrow 
{\rm{Exp}}_0tX$ which emanates from $o$ with initial direction $X$ in 
${\got m}_c$ is given by ${\rm{Exp}}_0tX=
\pi\exp tX$ , where $\exp$ is the exponential map from the Lie algebra
 ${\got g}_c$
to the Lie group $G_c$. A similar construction works in the case of $X_n$.

Below we specify the quantities introduced in this Section in the case of
 \Gras .
$o$ is taken as the $n$-plane {\bf O} given by eq. (\ref{z0}). The groups
are $G_c=SU(n+m),~ G_n=SU(n,m),~ G=SL(n+m,\gc ),~ K=S(U(n)\times U(m))$. The
symmetry at $o$ is $s= I_{nm}$, where $ I_{nm}=I_{nm}(-1)$ and
\begin{equation}
I_{nm}(\epsilon )=\left(\begin{array}{cc} \epsilon{\oo}_n & {\bf 0}\\
 {\bf 0} & {\oo}_m \end{array}\right)~,
\end{equation}
where $\epsilon =1~ (\epsilon=-1)$ for $X_c$ (resp. $X_n$). In fact
\begin{equation}
U^+I_{nm}(\epsilon)U=I_{nm}(\epsilon ),
\label{epse}
\end{equation}
where $U\in G_c $ ($G_n$) for $\epsilon=1$ (resp. $\epsilon =-1$).
We have also

\begin{equation}
{\got k}=\left\{\left ( \begin{array}{ll}
a & {\bf 0}\\
{\bf 0} & d
\end{array}
\right )|a^+=-a, d^+=-d, {\rm Tr}(a)+{\rm Tr}(d)=0\right\},
\end{equation}

\begin{equation}
{\got m}_{c,n}=\left\{\left (\begin{array}{ll}
{\bf o} & {\bf b} \\
-\epsilon b^+ & {\bf 0} 
\end{array}\right )\right\},
\end{equation}
 where $a$, $d$ and $b$ are, respectively, $n\times n$, $m\times m$ and
  $n\times m$ matrices.

The complex structure of the Grassmann manifold is inherited from his
representation as flag manifold\cite{wolf} $X_c=G^{\gcm}/P$,
 the parabolic group $P$ being 
\begin{equation}
P=\left\{ \left( \begin{array}{ll}
A & {\bf 0} \\ C & D
\end{array}
\right)
|~ \det A \det D =1
\right\}
\end{equation}
and $A$ $(D)$ is an $n\times n$ (resp. $m\times m$) matrix.

The compact roots are 
$$
\Delta_k= \left\{ e_i-e_j~|~1<i\not= j\leq n,~{\mbox {\rm or }} n<i\not= j
\leq m+n \right\}~,
$$
where $e_i,~i=1,\ldots ,N$ belong to the Cartan-Weyl basis.

The manifold $X_c$  and his noncompact dual $X_n$ can be parametrized as 
\begin{mathletters}\label{pxe}
\begin{eqnarray}
X_{n,c} & = & \exp\left(\matrix{{\bf 0}&B\cr
                           -\epsilon B^+&{\bf 0}\cr}\right) o=
   \left(\matrix{{\rm co}\sqrt{BB^+}&B{\displaystyle 
{{\rm si} \sqrt{B^+B}\over \sqrt{B^+B}}}\cr
                            -\epsilon{\displaystyle {{\rm si} \sqrt{B^+B}\over
                         \sqrt{B^+B}}}B^+&{\rm co}\sqrt{B^+B}\cr}\right)o
\label{57.1} \\*[2.0ex]
 & = & \left(\matrix{{\oo}&Z\cr {\bf 0}&{\oo}\cr}\right)
\left(\matrix{({\oo}+\epsilon ZZ^+)^{1/2}&
{\bf 0}\cr {\bf 0}&({\oo}+\epsilon Z^+Z)^{1/2}\cr}\right)
\left(\matrix{{\oo}&{\bf 0}\cr -\epsilon Z^+&{\oo}\cr}\right)o
\label{56}\\*[1.0ex]
\nonumber &  & \\
 & = & \exp\left(\matrix{{\bf 0}&Z\cr {\bf 0}&{\bf 0}\cr}\right)P ,\label{57}
\end{eqnarray}
\end{mathletters}
where {\rm co}  is an abbreviation for the circular cosine  cos 
(resp. the hyperbolic cosine coh) for $X_c$ (resp. $X_n$) and similarly for
 si.
 The sign  $\epsilon =+~(-)$  in eqs. (\ref{57.1}), (\ref{56}) corresponds to
 the compact (resp. noncompact) $X$.
In eq. (\ref{57}) $Z$ is the  $n\times m$ matrix of Pontrjagin
coordinates in ${\cal V}_0$ related to $B$ by the formula
\begin{equation}
Z=Z(B)=B{{\rm ta} \sqrt{B^+B}\over \sqrt{B^+B}} ,\label{58}
\end{equation}
and ta is an abbreviation for the hyperbolic tangent tanh (resp. the
circular tangent tan) for $X_n$ (resp. $X_c$) and eq. (\ref{58})
 realises the exponential
map in ${\cal V}_0$.

The noncompact case is realised under the restriction
\begin{equation}
{\oo}_m-Z^+Z>0~.\label{59}
\end{equation}

The representation (\ref{57})  for the noncompact case is the
 Harish-Chandra embed\-ding\cite{knapp}
of the noncompact dual $X_n$ of $X_c$   in $X_c$. Note that
 because of (\ref{57})
the complex matrix $Z$ parametrizes the Grassmann manifold.

The  invariant metric on \Gras ~to the group action, firstly studied by
 Teleman\cite{tel} and
Leichtweiss\cite{leich}, in the Pontrjagin coordinates reads
\begin{equation}
ds^2=k{\rm Tr}[(\oo +\epsilon ZZ^+)^{-1}dZ(\oo +\epsilon Z^+Z)^{-1}dZ^+].
\label{dist}
\end{equation}

$k=1$ in eq. (\ref{dist}) for $X_c=\Gras$ corresponds to
\begin{equation}
ds^2=\iota_*ds^2|_{FS}
\end{equation}
where
\begin{equation}
ds^2|_{FS}=d^2_c([\omega ],[\omega +d\omega ])=\frac{(\omega,\omega )(d\omega
 ,d\omega )-
(\omega ,d\omega )(d\omega , \omega )}{(\omega ,\omega )^2},
\end{equation}
and similarly for $X_n$.

The equation of the geodesics for $X_{c,n}$ is
\begin{equation}
\frac{d^2Z}{dt^2}-2\epsilon \frac{dZ}{dt}Z^+
(\oo +\epsilon ZZ^+)^{-1}\frac{dZ}{dt}=0~,\label{65}
\end{equation}
where $\epsilon =1~(-1)$ for $X_c$ (resp. $X_n$). It is easy to see that
$Z=Z(tB)$ in (\ref{58})
verifies (\ref{65}) with the initial condition $\dot{Z}(0)=B$.

A realization  of the algebra {\got a} consists of vectors of the form
\begin{equation}
H=\sum_{i=1}^rh_iD_{i\;n+i},~ h_i\in{\gr},\label{has}
\end{equation}
where $r$ is the symmetric rank of $X_c$ (and $X_n$) and we use the notation
\begin{equation}
D_{ij}=E_{ij}-E_{ji},~i,j=1,\ldots ,N.
\end{equation}
$E_{ij}$ is the matrix with entry $1$ on the $i$-th line and  $j$-th
column and $0$ otherwise.
We shall also need the notation
\begin{equation}
S_{ij}=E_{ij}+E_{ji}.
\end{equation}
 
The following Lemma will be used in order to calculate the tangent conjugate
 locus:
\begin{lem}
\label{lemm1}
The restricted roots of ($G_c, K$) are given in Table I, while the root
 space vectors are presented below.

\begin{table}
\caption[Table I]{{\rm Restricted roots for} \Gras . The lower indices of the
 roots vectors $X^j_{ab}$ are:  for 
$1\leq j\leq 8:~1\leq a\not=b\leq r;$ for $9\leq j\leq 12$: for every fixed
$a=1,\ldots ,r,~ b=1,\ldots ,|m-n| $ (if $m\not= n$); for $1\leq j\leq 14:~ a=
1,\ldots ,r .~( i=\sqrt{-1})$.}
\begin{center}
\begin{tabular}{||c|c|c||}\hline
{\rm Root space vectors} & {\rm Roots}  &{\rm Multiplicity} \\ \hline
$X^1_{ab};X^5_{ab}$ & $i(h_a-h_b)$ & $2 $ \\  \hline
$X^2_{ab};X^6_{ab}$ & $i(h_b-h_a)$ & $2 $ \\ \hline
$X^3_{ab};X^7_{ab}$ & $-i(h_a+h_b)$ & $2 $ \\ \hline
$X^4_{ab};X^8_{ab}$ & $i(h_a+h_b)$ & $2 $ \\ \hline
$X^9_{ab};X^{11}_{ab}$ & $ih_a$ & $2|m-n| $ \\ \hline
$X^{10}_{ab};X^{12}_{ab}$ & $-ih_a$ & $2|m-n| $ \\  \hline
$X^{13}_{a}$ & $2ih_a$ & $1 $ \\  \hline
$X^{14}_{a}$ & $-2ih_a$ & $1 $ \\  \hline
\end{tabular}
\end{center}
\end{table}

The first eight eigenvectors $X^{1-8}$ correspond to the eigenvalues
\begin{equation}
\lambda_{ab}=\epsilon _2i(h_a+\epsilon_1h_b),~ \epsilon_1^2=\epsilon_2^2=1,
\end{equation}
of the equation
\begin{equation}
\label{root}
[H,X^j_{ab}]=\lambda_{ab}X^j_{ab},~ \forall H\in{\got a},~ X\in{\got g}^{\gcm}.
\end{equation}

With the notation
\begin{equation}
X^j_{ab}=D^{\epsilon_1\epsilon_2}_{ab},~~
X^{j+4}_{ab}=S^{\epsilon_1\epsilon_2}_{ab},
\end{equation}
\begin{equation}
j=j(\epsilon_1 ,\epsilon_2)=\epsilon_1(1+\epsilon_2/2)+5/2,
\end{equation}
and if $F$ is any of the matrices $D$ and $S$, then
\begin{equation}
F^{\epsilon_1\epsilon_2}_{ab}=F_{a\,n+b}+
\epsilon_1F_{n+a\,b}+i\epsilon_2(F_{n+a\,n+b}-\epsilon_1F_{a\,b}).
\end{equation}

Explicitly
\begin{eqnarray}
\nonumber X^{1;5}_{ab} & = &~ i(F_{ab}+F_{n+a\,n+b})-F_{n+a\,b}+F_{a\,n+b},\\
\nonumber X^{2;6}_{ab} & = &-i(F_{ab}+F_{n+a\,n+b})-F_{n+a\,b}+F_{a\,n+b},\\
\nonumber X^{3;7}_{ab} & = &~i(F_{ab}-F_{n+a\,n+b})+F_{n+a\,b}+F_{a\,n+b},\\
\nonumber X^{4;8}_{ab} & = &- i(F_{ab}-F_{n+a\,n+b})+F_{n+a\,b}+F_{a\,n+b},
\end{eqnarray}
where the first (second) upper index of $X$ corresponds to
$F=D$ (resp. $F=S$).

The other vectors are as follows
$$X^{13}_a=\frac{1}{2}X^8_{aa},~X^{14}_a=\frac{1}{2}X^7_{aa},$$

$$~ X^{9;10}_{ab}  =   \left\{\begin{array}{ll}
E_{n+a\,2n+b}\mp iE_{2n+b\,n+a} & {\mbox{\rm if~}}n\leq m ,\\
E_{n+a\,m+b}\mp iE_{a\,m+b} & {\mbox{\rm if~}}n> m; 
\end{array}\right. $$

$$ X^{11;12}_{ab}  =  \left\{\begin{array}{ll}
E_{2n+b\,a}\pm iE_{2n+b\,n+a} &~~~{\mbox{\rm if~}}n\leq m, \\
E_{m+b\,a}\mp iE_{m+b\,n+a} &~~~{\mbox{\rm if~}}n> m.
\end{array}\right.$$

\end{lem}

{\it Proof\/}:~ The simplest proof is to solve the eigenvalue equation
 (\ref{root}).~\pre

6.2. The manifold ${\widetilde {\bf M}}_{n,c}$ of coherent
 states\cite{klauder} (in the sense of Perelomov\cite{per})
corresponding
to $X_{n,c}$ is introduced  in the notation of Ref. \onlinecite{sbl}.
 The  manifold of coherent
vectors is  the holomorphic line bundle associated to the character of the
parabolic subgroup $P$, with base the manifold of coherent states taken a
 homogeneous K\"ahlerian  manifold. The coherent
states are paramerized by a matrix $Z$ in front of the noncompact positive 
roots which appear at the exponent.\cite{sbl}

We shall prove the following

\begin{rem}
\label{rem7} The coherent vector $\vert Z,j_0>=\vert Z>$, 
 where $Z$ it is an $n\times m$  matrix,
 corresponds to the $n-$plane of ${\gc}^{n+m}$   parametrized by
 the Pontrjagin coordinates $Z$   in  ${\cal V}_0$  leading to
  $\hat{Z}= ({\oo}_nZ)$.

 Moreover, we have the equality of the scalar product of coherent vectors
 $<\cdot |\cdot >$  and the hermitian scalar product $((\cdot ,\cdot))$
 of the holomorphic line bundle $Det^*$:

\begin{equation}
<Z',j_0\vert Z,j_0>=<Z'\vert Z>=((Z',Z))=((Z^{\perp},Z'^{\perp}))~.\label{scal}
\end{equation}
 and similarly for the noncompact manifold $X_n$.
\end{rem}

{\it Proof}:  The scalar product of two coherent vectors  is\cite{sbl} 

\begin{equation}
<Z',j|Z,j>=\prod _{k=1}^{m+n}(A^{k\: k+1\ldots m+n}_{k\: 
k+1\ldots m+n})^{j_k-j_{k+1}}~,~j_0=0,~ j_1\geq j_2 \geq \cdots \geq j_{m+n} ~,
\end{equation}
where $A$ is the matrix

\begin{equation}
A= \left( \begin{array}{cc}
({\oo}_n+\epsilon ZZ'^+)^{-1} & {\bf 0} \\
{\bf 0} & ({\oo}_m+\epsilon Z'^+Z) 
\end{array}
\right) ~,
\end{equation}
the sign $\epsilon =-~ (+)$ corresponds to $X_n$ (resp. $X_c$) and the
 coherent vectors are considered in the chart ${\cal V}_0$ in the case of
 $X_c$.

Using the particular dominant weight

\begin{equation}
j=j_0=(\underbrace{1,\ldots ,1}_n,\underbrace{0,\ldots ,0}_m)~,
\end{equation}
it is found\cite{sbl} that
\begin{equation}
<Z',j_0\vert Z,j_0>=\det (\oo _n+\epsilon ZZ'^+)^{\epsilon } ~,\label{63}
\end{equation}
under the condition (\ref{scb}) of the scalar product.~\pre

Note that in the convention of Ref. \onlinecite{sbl}
 the coherent vector $|Z>$ corresponds to the $n-$plane $Z^t$
as a consequence of the fact that the fixed base point of the Grassmann 
manifold
 in Ref. \onlinecite{sbl} was chosen 
$$Z_0=\mb{e}_{m+1}\wedge\ldots\wedge\mb{e}_{n+m}$$ and not $Z_0$ given by 
(\ref{z0}). So, using Remark \ref{rem3}, it follows that eq. (\ref{scal})
should corresponds in the conventions of Ref.\onlinecite{sbl} to

\begin{equation}
<Z',j_0\vert Z,j_0>=<Z'\vert Z>=((Z'^{t},Z^{t}))=((Z^{t\perp},Z'^{t\perp}))~.
\end{equation}

Finally, we remember that Calabi's diastasis function\cite{cal} $D(Z',Z)$ has
been used in the  context of coherent states,\cite{cgr} observing that 
$D(Z',Z)=-2\log
<\underline{Z'}|\underline{Z}>$, where $|\underline{Z}>=<Z|Z>^{-1/2}|Z>$.

The noncompact Grassmann manifold $X_n$ admits the embedding  in an 
infinite
dimensional projective space $\iota_n:X_n\hookrightarrow {\bf PK}$ and
 also the embedding
 $\iota ':X_n\hookrightarrow {\gcp}^{N(n)-1,1}$. Let $\delta_n$
($\theta_n$) be the length of the geodesic joining $\iota '(Z'),\iota '(Z)$
 (resp. $\iota_n(Z'),\newline\iota_n(Z)$).  Then we have the 
\begin{rem}
\label{remc}
For the noncompact Grassmann manifold, $\delta_n,~\theta_n$ and $D_n$
are related through the relation
\begin{equation}
\cos\theta_n=\cosh^{-1}\delta_n=e^{-D_n/2}=\frac{|\det (\oo -ZZ'^+)|}{\det
[(\oo -ZZ^+)(\oo -Z'Z'^+)]^{1/2}}~.
\end{equation}
\end{rem}
{\it Proof}:  
Due to eq. (\ref{scal}), the diastasis is related to the geodesic distance
$\theta (Z',Z) $ given by eq. (\ref{10.4.1}) by the relation $D(Z',Z)=
-2\log \cos (\theta (Z',Z))$. So, the diastasis for $X_c$ if $Z,
 Z'\in{\cal{V}}_0$ is
\begin{equation}
D(Z',Z)=\log\frac{\det[ (\oo +ZZ^+)(\oo +Z'Z'^+)]}{|\det (\oo +ZZ'^+)|^2},
 \label{cah}
\end{equation}
and similarly for the noncompact case.~\pre

The equation (\ref{cah}) is still valid for the infinite dimensional Grassmann
manifold.\cite{ms}

\section{\hspace{-4mm}.\hspace{2mm}  CUT LOCUS AND CONJUGATE  LOCUS}
\hspace*{0.8cm} 
\subsection{\hspace{-5mm}.\hspace{2mm} Preliminaries}

 We begin remembering some definitions referring to the cut
locus and conjugate locus.

 Let $V$ be a compact Riemannian manifold of dimension $n$, $p\in V$ and let
 ${\rm Exp}_p$
be the (geodesic) exponential map at the point $p$. Let $C_p$ denote
the set of vectors $X\in V_p$ (the tangent space at $p\in V$) for which
${\rm Exp}_pX$ is singular. A point $q$ in $V$ ($V_p$) {\it is conjugate
to p} if it is in ${\bf C}_p={\rm Exp}\,C_p$ ($C_p$)\cite{helg} and ${\bf
C}_p$ ($C_p$) is called {\it the conjugate locus} (resp. {\it tangent
 conjugate locus}) of the point $p$. 

Let $q\in V$. The point $q$ is in the {\it cut locus} ${\bf CL}_p$ of $p\in V $
if it is nearest point to $p\in V$ on the geodesic joining $p$ with $q$,
beyond which the geodesic ceases to minimize its arc length.\cite{kn}  More
 precisely, let $\gamma _X(t)={\rm Exp}\, tX$ be a geodesic emanating from
 $\gamma_X(0)=p\in V$, where  $X$ is a unit vector from the unit sphere $S_p$
 in $V_p$.
$t_0X$ (resp. ${\rm Exp}\,t_0X$) is called a {\it tangential cut point}
 ({\it cut point}) of $p$ along $t\rightarrow{\rm Exp}\,tX$
 ($0\leq t\leq s$) if
the geodesic segment joining $\gamma_X(0)$ and $\gamma_X(t)$ is a minimal
geodesic for any $s\leq t_0$ but not for any $s>t_0$.

 Let us define the function $\mu :
S_p\rightarrow {\gr}^+\cup \infty,~ \mu (X)=r$, if $q={\rm Exp}\,rX\in
{\bf CL}_p$, and $\mu (X)=\infty$ if there is no cut point of $p$ along 
$\gamma_X(t)$. Setting $ I_p=\{tX,~0<t<\mu (X)\}$, then ${\bf I}_p=
{\rm Exp}\,I_p$ is called the {\it  interior set at p}. Then\cite{kn}

1) ${\bf I}_p\cap {\bf CL}_p=\emptyset ,~ V={\bf I}_p\cup {\bf CL}_p,$ the
 closure $\bar {\bf I}_p=V$, and dim ${\bf CL}_p\leq n-1$.

2) $I_p$ is a maximal domain  containing $0=0_p\in V_p$ on which
 ${\rm Exp}_p$ is a
 diffeomorphism and ${\bf I}_p$ is the largest open subset of $V$ on which a 
normal coordinate system around $p$ can be defined.

This theorem will be used below in the proof of Proposition \ref{Prop3}.

The importance of the cut loci lies in the fact they inherit topological
 properties of the manifold $V$.

The relative position of ${\bf CL}_0$ and ${\bf C}_0$ is given by Theorem 7.1
 p. 97 in Ref. \onlinecite{kn} :

 Let the notation
 $\gamma_t=\gamma _X(t)$. 
Let $\gamma_r$ be the cut point of $\gamma_0$ along a geodesic $\gamma = 
\gamma_t,~ 0\leq
t< \infty$. Then, at least one (possibly both) of the following
statements holds:

(1) $\gamma_r$ is the first conjugate point of $\gamma_0$ along $\gamma$;

(2) there exists, at least, two minimising geodesics from $\gamma_0$ to
 $\gamma_r$.

Crittenden\cite{cr} has shown that for the case of simply connected symmetric
spaces, the cut locus is identified with the first conjugate locus. This result
will be illustrated on the case of the complex Grassmann manifold. Generally,
 the situation  is more complicated.\cite{war,wei}

 For ${\gcp}^n$, $CL$ is
the sphere of radius $\pi$ with centre at the origin of the tangent space to
${\gcp}^n$ at the given point, while ${\bf CL}$ is the hyperplane at infinity
${\gcp}^{n-1}$. Except few situations, e. g.  the ellipsoid, even for low
 dimensional manifolds, ${\bf CL}$
is not known explicitly. Helgason\cite{helg} has shown that the  cut locus of
 a
compact connected Lie group, endowed with a bi-invariant Riemannian metric is 
stratified, i.e. it is the disjoint union of smooth submanifolds of $V$. This
situation will be illustrated on the case of the complex Grassmann manifold.
Using a geometrical method, based on the Jordan's stationary angles,
 Wong\cite{wong,won,wo} has studied conjugate loci
 and cut loci of the Grassmann manifolds.  Calculating the tangent 
conjugate locus on the Grassmann manifold, Sakai\cite{sak}
observed that the results of Wong\cite{won} referring to conjugate locus in
Grassmann  manifold are incomplete. This problem will be largely discussed in
 the present Section. By refining the results of Ch. VII,
 \S 5 from Helgason's book,\cite{helg} Sakai\cite{sak,sak1}
studied the  cut locus on a general symmetric space and showed that it is
determined by the cut locus of a maximal totally geodesic flat submanifold of
 $V$. However, the expression of the conjugate locus as subset of the Grassmann
manifold is
not known explicitly. We give a geometric characterization of the part of the
conjugate locus different from those found by Wong in terms of the stationary
angles.
\subsection{\hspace{-5mm}.\hspace{2mm} Cut locus} 
Coming back to eqs. (\ref{57.1})-(\ref{57}), it is observed that
 $B$ {\it are normal
coordinates around} $Z=0$ on the Grassmann manifold. So we have
\begin{equation}
<{\bf O}|Y>=0 ~{\mbox {\rm iff }}~ (({\bf O},Y))=0~ {\mbox {\rm or iff}}~ Y\in
\Sigma _0 ~.
\end{equation}

The following two assertions are particular situations true for symmetric
or generalized symmetric spaces\cite{sb,sb1}

\begin{pr}[Wong\cite{wong}]
\label{Prop3}  The cut locus, the polar divisor of ${\bf O} \in 
{\cal V}_0 \subset G_n({\bf K})$ and the interior set  are related by the
 relations
\begin{equation}
{\bf CL}_0=\Sigma_0~,
\end{equation}
\begin{equation}
{\cal V}_0={\bf I}_0,
\end{equation}
and $\Sigma_0$ is given by Lemma \ref{lemma3}.
\end{pr}

\begin{rem}
\label{remark9}
 The solution of the equation
\begin{equation}
<0|\psi >=0~,
\end{equation}
 where $|\psi >$ is a coherent vector, is given by the points on the
 Grassmann manifold corresponding to the cut locus  ${\bf CL}_0=\Sigma_0$.

\end{rem}

{\it Proof}: The dependence $Z(t)=Z(tB)$ expressed by  (\ref{58})
  gives geodesics starting at
 $Z=0$ in the
chart ${\cal V}_0$ and ${\cal V}_0$ {\it is the maximal normal neighbourhood}.
The Proposition follows due to Thm. 7.4 of Kobayashi and Nomizu\cite{kn} and
 the subsequent remark  at p. 102 reproduced earlier.~\pre

\subsection{\hspace{-5mm}.\hspace{2mm} The conjugate locus in the complex
 Grassmann manifold} 

 Now the conjugate locus of the point $Z=0$ in \Gras ~ is calculated. The 
Jacobian
of the transformation (\ref{58}) has to be computed. We shall prove the
following theorem and remark

\begin{thm}
\label{ws}
 The conjugate locus of {\bf O} in \Gras  ~is given by the union
\begin{equation}
\label{reun}
{\bf C}_0={\bf C}^W_0\cup {\bf C}^I_0.
\end{equation}

${\bf C}^W_0$ consists of those points of the Grassmann manifold which have at
least one of the stationary angles with the {\bf O} plane $0$ or $\pi /2$.
${\bf C}^I_0$ consists of those points of \Gras  \ for which at least two of
 the stationary angles with {\bf O} are equal, that is  at least two of the
 eigenvalues of the matrix (\ref{53}) are equal.

The ${\bf C}^W_0$ part of the conjugate locus is given by the disjoint union
\begin{equation}
{\bf C}^W_0=\label{72}
\cases {V^m_1\cup V^n_1,&$ n\leq m,$\cr
            V^m_1\cup V^n_{n-m+1},& $n>m,$ \cr}
\end{equation}
 where
\begin{equation}
V^m_1=\label{74.1}
\cases {{\gcp}^{m-1},&  $n=1 ,$\cr
              W^m_1\cup W^m_2\cup \ldots W^m_{r-1}\cup W^m_r,&$1<n ,$\cr}
\end{equation}
\begin{equation}
W^m_r=\cases {G_r({\gc}^{\max (m,n)}),&$n\not= m,$\cr
             {\bf O}^{\perp},&$n=m ,$\cr}
\end{equation}
\begin{equation}
V^n_1=\label{72.2}
\cases {W^n_1\cup \ldots \cup W^n_{r-1}\cup {\bf O},&$1<n\leq m,$\cr
             {\bf O}, &$n=1 ,$ \cr}
\end{equation}
\begin{equation}
V^n_{n-m+1}=W^n_{n-m+1}\cup W^n_{n-m+2}\cup\ldots\cup W^n_{n-1}\cup
{\bf O}~,~ n>m ~.\label{73.3}
\end{equation}
\end{thm}

\begin{rem}
\label{remr10}
The cut locus in \Gras \ is given by those $n$-planes  which have at least one
angle $\pi /2$ with the plane {\bf O}.
\end{rem}

{\it Proof\/}:~ The proof is done in 4 steps. $a)$ Firstly, a diagonalization 
is performed.
 $b)$ After this, the Jacobian of a transformation  of complex dimension one
 is computed. At $c)$
 the cut locus is reobtained. d) Finally, the property of the stationary
 angles given in Comment 1 is used in order to get the conjugate locus in
 \Gras .

At $b)$ we argue that the proceeding gives all the conjugate locus. See also the
proof of the Proposition \ref{sakth}, where it is stressed
 the equivalence of the
decomposition (\ref{hu}), (\ref{68}) with the representation given by eq.
 (\ref{helgas}).

$a)$ Every $n\times m$ matrix $Y$ can be put in the form\cite{hua}
\begin{equation}
Y=U\Lambda V,\label{hu}
\end{equation}
 where $U~(V)$ is a unitary $n\times n$ (resp. $m\times m$) matrix,
\begin{equation}
\Lambda=\label{68}
\left\{ \begin{array}{c@{\quad {\mbox {\rm if~}}\quad}l}
(D,0) & n\leq m~, \\*[1.0ex]
\left( \begin{array}{c}
\! D \! \\ \! 0 \! 
\end{array}
\right)
 & n>m~,
\end{array}
 \right.
\end{equation}
and $D$ is the $r\times r$ diagonal matrix $(r=\min (m,n))$, with diagonal
 elements $\lambda _i \geq 0, ~ i=1,\ldots , r$. If the rank of the matrix $Y$
is $r_1$, then $\Lambda$ has four block form with the diagonal elements
 $\lambda_i>0,~i=1,\ldots ,r_1$, the other elements being $0$.\cite{ben}

We shall apply a decomposition of the type (\ref{hu})-(\ref{68}) taking the
 diagonal
elements of the matrix  $D$ as complex numbers. This implies an overall phase
$e^{i\varphi}$ for the matrix $\Lambda$ in the decomposition (\ref{hu}).
 This phase can be included in the matrix $U$ such that the matrix $U'=
e^{-i\varphi}U$ is unitary in the decomposition (\ref{hu}).

Applying to the $n\times m$ matrix $B$ the diagonalization technique here
 presented, the $n\times m$ matrix  $Z=Z(tB)$
 corresponding to eq. (\ref{58}) is also of the same form, where the
 diagonal elements  are 
\begin{equation}
Z_i=\frac{B_i}{|B_i|}\tan t|B_i|, ~i=1,\ldots , r~.\label{70}
\end{equation}

$b)$ In order to calculate the Jacobian $J$ when $B$ is diagonalized,
 let us firstly calculate
\begin{equation}
\Delta _1=\Delta _1(X,Y)= \label{71.1}\left| \begin{array}{ll}
{\displaystyle \frac {\partial X}{\partial B_x}} & 
{\displaystyle \frac {\partial X}{\partial B_y}} \\*[2.0ex]
{\displaystyle \frac {\partial Y}{\partial B_x}}   &
{\displaystyle  \frac {\partial Y}{\partial B_y}}
\end{array}
\right|,
\end{equation}
where $Z=X+iY,~ B=B_x+iB_y,$ and $ X,Y,B_x,B_y \in {\gr}$.

With eq. (\ref{70}) we get
\begin{equation}
\Delta _1={\displaystyle \frac {t}{|B|}\frac{\sin t|B|}{\cos ^3 t|B|}}~.
\end{equation}

 Now there are two possibilities.
 $\alpha )$ Firstly, we consider the case when all the $|B_i|$ in eq.
  (\ref{70}) are
 distinct. Then the Jacobian $J$ corresponding to the transformation (\ref{hu}),
 (\ref{68}) of $Z$ is \\

\begin{equation}
\Delta =\prod _{i=1}^r{\displaystyle \frac {t}{|B_i|}\frac{\sin t|B_i|}
{\cos ^3 t|B_i|}}~.\label{73}
\end{equation}

$\beta )$ Otherwise, at least two of the eigenvalues $|B_i|$
 in the transformation ({\ref{70}) are equal. Then the Jacobian $J$ is zero and
 the points are in the tangent conjugate locus.

$c)$ Let us now take in eq. (\ref{23.2}) (and (\ref{63})) $Z'=0$ and
 let us introduce $Z$ in the  diagonal form (\ref{68}), (\ref{70})
 in eqs. (\ref{53}), (\ref{54}) and 
(\ref{55}). It results  that 
\begin{equation}
\label{jr}
\cos\theta_i=|\cos\,t|B_i||,
~i=1,\ldots ,r ~,
\end{equation}
 and then formula (\ref{55})  becomes
\begin{equation}
\cos \theta = \prod_{i=1}^r|\cos t|B_i|| ~.\label{75}
\end{equation}

If for some $i\not= j$, $|B_i|$ and $|B_j|$ are identical or
\begin{equation}\label{peste}
|B_i\pm B_j|=\pi,
\end{equation}
then they correspond to the same stationary angles $\theta_i=\theta_j$, cf. eq.
(\ref{jr}). In other words, if at least two of the eigenvalues of the matrix 
(\ref{53}) are identical, then they correspond to the same stationary angles 
$\theta_i=\theta_j$.

If for some $i,~ t|B_i|=\pi /2$ in the  $Z$  matrix put in the diagonalized 
form,
we have to change the chart. As a consequence of the fact that the change of
 charts has the homographic form (\ref{3xxx}), a change of those coordinates
which  are not finite in one chart has the form $Z\rightarrow 1/Z$, the matrix
$B$ being diagonalized.
 So, we have to calculate instead of eq.
 (\ref{71.1}),
the Jacobian $\Delta '_1=\Delta '_1(X',Y')$, where $Z'=1/Z=X'+iY', X',~Y'\in
\gr$.
It is easily found out that
\begin{equation}
\Delta '_1={\displaystyle \frac {t}{|B|}\frac{\cos t|B|}{\sin^3 t|B|}}~.
\label{71.2}
\end{equation}

With equation (\ref{71.2}) we get that $\Delta'=0$ iff
 $\cos \theta =0$ in eq. (\ref{75}),
where $\Delta '$ is the expression corresponding to $\Delta$ in the new
chart. In fact, $\cos \theta =0$ iff for at least  one $i$,  $t|B_i|=\pi /2,
~i=1,\ldots ,r$, or, equivalently, iff at least one of the angles of
 ${\bf O}^\perp$
 with $Z$ is zero, i.e. $\dim ({\bf O}^\perp\cap Z)\geq 1$. The
results  of Proposition \ref{Prop3} and Remark \ref{remark9}  for the cut
 locus are reobtained,
 i.e. ${\bf CL}_0= \Sigma_0=V^m_1$ and the Remark \ref{remr10} is proved.

From formulas (\ref{73}), (\ref{75}) it follows also that {\it in the tangent
 space the cut locus }= {\it  first
conjugate locus}, a result true for any symmetric simply connected
 space\cite{cr} as has been already remarked.

$d)$  Further we look for the other points $Z$ in  the conjugate locus ${\bf
C}_0$ 
but not in ${\bf CL}_0$, i.e. we look for the other points where $J=0$. 

Once the cut locus was gone beyond,
  the same chart as before the cut locus has
been reached can be used. Then at lest one of the $t|B_i|$ is zero 
(modulo $\pi$),
corresponding to at  least
one of the angles between {\bf O} and $Z$  zero.

Let $n\leq m$.  If $Z\in G_n({\gc}^{n+m})$ is such that 
$\dim (Z\cap {\bf O})=i, ~ i=1,\ldots n,$ then $i$  stationary angles
between {\bf O} and $Z$  are zero ($n-i$
are different from 0) and $Z\in W^n_i$. Finally $Z\in {\bf C}_0\setminus 
{\bf CL}_0$
iff
 $$Z\in \bigcup_{i=1}^nW^n_i=V^n_1 ~.$$

Let now $n>m$. We look for points for which $J=0$ different from the points of
${\bf CL}_0$ where $\Delta=0$. So eq. (\ref{73}) with $r=m$ does not include 
the $n-m$ stationary angles
which are already zero. Then $J=0$ when at least $n-m+1$ angles are zero. In
fact, if {\bf O} and $Z$ are two $n-$planes such that $n-m+i$ angles are
 zero (and $m-i$ angles are not zero), then $\dim ({\bf O}\cap Z)=n-m+i, ~ i=1,
\ldots ,m$. So, in order that the equation $J=0$ to be satisfied  in the
 points of ${\bf C}_0\setminus {\bf CL}_0$, it is necessary and sufficient that
$$Z \in \bigcup _{i=1}^m W^n_{n-m+i} = V^n_{n-m+1}~. $$

The representations (\ref{74.1}), (\ref{72.2}), (\ref{73.3})  follow
 particularizing  the third eq. (\ref{43})
and the last term is obtained as particular case of eq. (\ref{44}).

To see that the union (\ref{72}) is disjoint, it is observed that
 $V^m_1=Z(m-1,m,\ldots ,m)$ and  $V^n_1=Z(n-1,m,\ldots ,m)$.
The condition to have nonvoid intersection of the Schubert varieties
$Z(\omega )~, Z(\omega ')$ is that $\omega_i+\omega'_{n-i}\geq n+m,~
i=1,\ldots ,n$
(cf.   p. 326 in Ref. \onlinecite{hodge}).~\pre

 Wong's\cite{won} notation is
\begin{equation}\label{wrin}
V_l=\{Z\in G_n({\bf C}^{n+m})  |
\dim (Z\cap {\bf O}^\perp )\geq l \},
~{\widetilde V}_l=\{Z\in G_n({\bf C}^{n+m})  |
\dim (Z\cap {\bf O}) \geq l \}~,
\end{equation}
i.e. $V_l~(\widetilde{V}_l)$ from Wong corresponds to our $V^m_l$
(resp. $V^n_l$) and
\begin{equation}
{\bf C}^W_0=
\left\{
\begin{array}{l@{\quad , \quad}l}
V_1\cup {\widetilde V}_1 & n\leq m ~, \\
V_1\cup {\widetilde V}_{n-m+1} & n> m ~.
\end{array}
\right. 
\end{equation}
\subsection{\hspace{-5mm}.\hspace{2mm} The tangent conjugate locus}

\indent
The tangent conjugate locus $C_0$  for \Gras~  in
the case $n\leq m$ was obtained by
 Sakai.\cite{sak} Sakai has observed that Wong's result on the conjugate locus
in the manifold is incomplete, i.e. ${\bf C}^W_0\subset{\bf C}_0$ but
${\bf C}^W_0\ine {\bf C}_0=\exp C_0$. The proof of Sakai consists in
solving the eigenvalue equation $R(X,Y^i)X=e_iY^i$
which appears when solving  the Jacobi equation,
 where the
curvature for the symmetric space $X_c=G_c/K$ at $o$ is simply $R(X,Y)Z=
[[X,Y],Z],~X,Y,Z\in\got{m}_c$. Then $q={\rm Exp}_0tX$ is conjugate to $o$ if
$t=\pi \lambda/\sqrt{e_i},~\lambda\in\gz^{\star}\equiv \gz\setminus \{0\}$.
 The solution of the same
problem in terms of $\alpha (H)$ is given by Lemma 2.9 at page 288 in
 Ref. \onlinecite{helg}. In \S 7.3 we have calculated ${\bf C}_0$ using
 directly the form (\ref{58}) of the exponential map in ${\cal V}_0$.
 Below we present another calculation of $C_0$ and compare
these  results with those proved in Theorem
 \ref{ws} referring to the conjugate locus in \Gras .

\begin{pr}
\label{sakth}
The tangent conjugate locus $C_0$ of the point ${\bf O}\in\Gras$ is given by 
\begin{equation}
\label{ura}
C_0=\bigcup_{k,p,q,i}ad\,k(t_iH)~,~i=1,2,3;~1\leq p<q\leq r,
~ k\in K,
\end{equation}
where the vector $H\in\got{a}$  (eq. \ref{has}) is normalized,
\begin{equation}
H=\sum_{i=1}^r h_iD_{i\,n+i},~h_i\in\gr,~\sum h^2_i=1~.
\label{hh}
\end{equation}
 The parameters $t_i,~i=1,2,3$ in eq. (\ref{ura}) are
\begin{equation}
\begin{array}{l@{\:=\:}c@{\:,\:}l}
t_1  & \displaystyle{\frac{\lambda \pi}{|h_p\pm h_q|}}  &
 ~\mbox{\rm multiplicity}~2;\\[2.ex]
t_2  & \displaystyle{\frac{\lambda \pi}{2|h_p|}} & 
  ~\mbox{\rm multiplicity}~1;\\[2.ex]
t_3  & \displaystyle{\frac{\lambda \pi}{|h_p|}}  & 
 ~\mbox{\rm multiplicity}~2|m-n|; ~\lambda\in \gz^{\star}~.
\end{array}
\label{ttt}
\end{equation}
 The following relations are true

\begin{equation}
\label{ect1}
{\bf C}^I_0= \exp \bigcup_{k,p,q}Ad\,k(t_1H)~,
\end{equation}
\begin{equation}
\label{ect2}
{\bf C}^W_0= \exp \bigcup_{k,p}Ad\,k(t_2H)~,
\end{equation}
i.e. exponentiating the vectors of the type $t_1H$ we get the points
of ${\bf C}^I_0$ for which at least two of the stationary angles with {\bf O}
are equal, while the vectors of the type $t_2H$ are sent to the points of
 ${\bf C}^W_0$ for which at least one of the stationary angles with {\bf O}
is $0$ or $\pi /2$.
\end{pr}

{\it Proof\/}:~ Any vector $X\in\got{m}_c$ can be put\cite{helg} in the form
\begin{equation}
\label{helgas}
X=Ad\,(k)H,~k\in K,~H\in\got{a}.
\end{equation}

So, in order to find out the tangent conjugate locus it is sufficient to solve
this  problem for $H\in\got{a}$. But then $X$ in eq. (\ref{helgas}) is
conjugate with $o$ iff
\begin{equation}
\label{alphac}
\alpha (H)\in i\pi \gz^{\star}
\end{equation}
for some root $\alpha$ which do not vanishes identically on \got{a} (cf.
i.e.  Prop. 3.1 p. 294 in the book of Helgason\cite{helg}). In fact, what we
 have to find out is the diagram of the pair $(G_c,K)$.

But according to Lemma \ref{lemm1}, the space of restricted root vectors
consists of three types of vectors, corresponding respectively to the
 restricted roots: $\pm i(h_a-h_b), \pm ih_a$ and $\pm 2ih_a$, where $1\leq a<
b\leq r$, with the vector $H$ of the form (\ref{hh}). So, imposing to the
 vectors $tH$ the condition (\ref{alphac}), the values
 (\ref{ttt}) are obtained for the parameters $t_i$.

To compare the results on $C_0$ with those on ${\bf C}_0$, let us observe
 that a diagonal matrix $B$ as in \S 7.3 corresponds to the representation
 (\ref{hh}). When expressed in stationary angles, the "singular value
 decomposition"  (\ref{hu}) is nothing else than the representation
 (\ref{helgas}) expressed matricially. In eq. (\ref{68})
 $\Lambda$ corresponds to $B$ while $D$ corresponds
to $H$, where $B$ is in $\got{m}_c$ as in eq. (\ref{57.1}) and $H$ has the form
 (\ref{hh}). This implies that for the vector
\begin{equation}
\label{v100}
t_1H=\frac{\pi \lambda }{|h_p\pm h_q|}\sum_{i=1}^r h_iD_{i\,n+i},~\lambda\in
\gz^{\star}~,
\end{equation}
the $p^{th}$ ($q^{th}$) coordinate in the tangent space $\got{m}_c$ to
\Gras ~ at  $~o$ is $B_p=\lambda \pi h_p|h_p\pm h_q|^{-1}$
 (resp. $B_q=\lambda \pi h_q|h_p\pm h_q|^{-1}$). Consequently,
 the relation (\ref{peste}) is fulfilled .
  So, due to eq. (\ref{jr}), the corresponding stationary
angles are equal, $\theta_p=\theta_q$ and eq. (\ref{ect1}) is proved.

Similarly, the vector
\begin{equation}
\label{v101}
t_2H=\frac{\pi \lambda }{2|h_p|}\sum_{i=1}^r h_iD_{i\,n+i},~\lambda\in
\gz^*~,
\end{equation}
corresponds for $\lambda$ even (odd) to points on \Gras \ which have at least 
one of the stationary angles with {\bf O} equal to $0$ (resp. $\pi/2$). This
fact and also the representation (\ref{72}) can be seen with eq. (\ref{57.1})
with diagonal $B$-matrix. Then in eq. (\ref{57.1})
\begin{mathletters}
\begin{equation}
{\rm diag}(\sqrt{BB^+})=(|h_1|,\ldots ,|h_n|),~{\rm if}~n\leq m~,
\end{equation}
\begin{equation}
{\rm diag}(\sqrt{B^+B})=(|h_1|,\ldots ,|h_m|),~{\rm if}~n> m~,
\end{equation}
\end{mathletters}
where ${\rm diag}(X)$ denotes the diagonal elements of the matrix $X$.

 Choosing $o\in\Gras$ to correspond to {\bf O}
given by (\ref{z0}), then (cf. eq. (\ref{15}))  a point of
 {\bf O} has the coordinates 
$$(x_1,\ldots ,x_n,\underbrace{0,\ldots ,0}_{m})~.$$
So, a point of the Grassmann manifold $X_c$ is
\begin{mathletters}\label{las}
\begin{equation}
\label{last}
(x_1\!\cos\! |h_1|,\ldots\! ,x_n\!\cos\! |h_n|,-x_1\!\frac{h_1}{|h_1|}\sin\!
 |h_1|,\ldots\! ,
-x_n\!\frac{h_n}{|h_n|}\sin\! |h_n|,\underbrace{0,\ldots\! ,0}_{m-n}),
\end{equation}
\begin{equation}
\label{last1}
(x_1\!\cos\!|h_1\!|,\ldots\!,x_m\!\cos\! |h_m\!|,x_{m+1}\!,\ldots\!,x_n,
\!-x_1\!\frac{h_1}{|h_1\!|}\sin\! |h_1\!|,
\ldots\! ,
\!-x_m\!\frac{h_m}{|h_m\!|}\sin\! |h_m\!|),
\end{equation}
\end{mathletters}
where eq. (\ref{last}) ((\ref{last1})) corresponds to the case
  $~n\leq m$ (resp. to $n> m$).~\pre

Note that ${\bf C}_0^I$ is not a Schubert variety, because in general
 ${\bf C}_0^I\cap{\bf O}=\emptyset$ and  ${\bf C}_0^I\cap{\bf O}^{\perp}
 =\emptyset$. However,
the representation (\ref{v100}) with $h_p+h_q$ at the denominator gives $V_2$
($\widetilde{V}_2$) for $\lambda$ odd (resp. even) cf. eq. (\ref{last}) or
 Sakai.\cite{sak} But $V_2\subset V_1$ (in fact $V_2$ is the singular locus of
$V_1$) and the union in eq. (\ref{ura}) is not disjoint even for the parts
 which correspond
to $t_1H$ and $t_2H$.
\begin{com}
\label{com2}
${\bf C_0}^I$ contains as subset the maximal set of mutually isoclinic
subspaces of the  Grassmann manifold, which are the isoclinic spheres, with
dimension given by the solution of the Hurwitz problem.
\end{com}

{\it Proof}: Wong\cite{wwong} has found out the locus of isoclines in
 $G_n({\gr}^{2n})$,
i.e. the maximal subset $B$   of the Grassmann manifold containing {\bf O}
 consisting of points with
the property that every two $n$-planes   of $B$
have {\it all} the stationary angles equal. Two mutually isoclinic $n$-planes
correspond to the situation where the matrix (\ref{51}) is a multiple of
 {\oo}.  The results of Wong were generalized by Wolf,\cite{ww1} who has
considered also the complex and quaternionic Grassmann manifolds.
 The problem of maximal mutually isoclinic subspaces is related with the
Hurwitz problem.\cite{hur} Any maximal set of mutually isoclinic $n$-planes is
 analytically  homeomorphic to a sphere (cf. Thm 8.1 in Wong\cite{wwong} and
 Wolf\cite{ww1}), the dimension of the isoclinic spheres being given by the 
solution to the Hurwitz problem.~\pre

\section{\hspace{-4mm}.\hspace{2mm}  THE DISTANCE }
\hspace*{0.8cm}

In this Chapter $Z$ is an $n\times m$ matrix characterizing a point in the
 complex
Grassmann manifold $X_c$ (\ref{xc}) (resp. the noncompact dual (\ref{xn})
$X_n$ of $X_c$). In the case of the compact Grassmann manifold
 $X_c$ 
the $n-$plane $Z$ is taken in ${\cal V}_0$, while in the case of the noncompact
manifold $X_n$, the matrix $Z$ is restricted by the condition (\ref{59}). In 
formulas below $\epsilon =1$ (-1) and arcta is an abbreviation for the 
inverse of the circular tangent function, arctan (hyperbolic function arctanh)
for $X_c$ (resp. $X_n$) and analogously for arcco and arcsi.

Let us denote by $\lambda _i(A), ~i=1,\ldots ,p$ the eigenvalues of the
$p\times p$ matrix $A$ and let $\eta =\sqrt{-\epsilon}$.

We shall prove the following 

\begin{pr}
\label{pr5} The square of the distance between two points
$Z_1,~Z_2$  is given by the formulas
\begin{eqnarray}
d^2(Z_1,Z_2) & = &\mbox {\rm Tr} \left[\mbox {\rm arcta}(ZZ^+)^{1/2}\right]^2,
\label{tr} \\
             & = & \sum_{j=1}^n \theta_j^2~, \label{r}
\end{eqnarray}
 where

\begin{eqnarray}
Z & = & (\oo +\epsilon Z_1Z^+_1)^{-1/2}(Z_2-Z_1)(\oo +\epsilon Z^+_1Z_2)^{-1}
(\oo +\epsilon Z^+_1Z_1)^{1/2},\label{zz} \\
\theta _j & = & {\rm arcta}\,\lambda _j(ZZ^+)^{1/2} 
 =   {\rm arcco}\,\lambda _j(V)^{1/2}
\label{un}              =   {\rm arcsi}\,\lambda _j(VZZ^+)^{1/2}\\
\label{lag}  & = & \frac {1}{2\eta}\log \frac {\oo +\eta\lambda _j(ZZ^+)^{1/2}}
{\oo -\eta\lambda _j(ZZ^+)^{1/2}}=\frac{1}{\eta}\log\lambda_j
[V^{1/2}(\oo +\eta (ZZ^+)^{1/2})],\\
V & \equiv & (\oo +\epsilon ZZ^+)^{-1}, \label{vv}\\
V & = & (\oo\!+\!\epsilon Z_1Z^+_1)^{-1/2}(\oo\!+\!\epsilon Z_1Z^+_2)
(\oo\!+\!\epsilon Z_2Z^+_2)^{-1}(\oo\!+\!\epsilon Z_2Z^+_1)
(\oo\!+\!\epsilon Z_1Z^+_1)^{-1/2}\!.\label{v}
\end{eqnarray}

 The matrices $V$  and $W$  given by eq. 
(\ref{53n})  have the same eigenvalues.
\end{pr}

{\it Proof\/}:~ The proof is done in three steps. {\it a)} Firstly, a
 homographic 
 transformation $ Z'=Z'(Z)$ for which $Z'(Z_1)=0$ is obtained. {\it b)} Further
on, the distance between $Z=0$ and $Z$, where $Z\in {\cal V}_0$ for $X_c$
is found. {\it c)} Finally, the transformation $Z'=Z'(Z_2)$ gives eq.
 (\ref{zz}), while (\ref{un})-(\ref{lag}) are obtained as $d^2(0,Z'(Z_2))$.
The representation (\ref{v}) of the matrix (\ref{vv}) is furnished by a matrix
 calculation.

{\it a)}   The transitive action of an element from the group $G_c=SU(n+m)~
(G_n=SU(n,m))$ on $X_c$ (resp. $X_n$) is given by the linear fractional
transformation
       
\begin{equation}
Z'=Z'(Z)=U\cdot Z=
(AZ+B)(CZ+D)^{-1},~U=\left(\begin{array}{cc}A & B \\ C & D\end{array}\right)\in
G_c~(G_n)~.
\label{omo}
\end{equation}

So, we have to find a matrix $U\in G_c ~(G_n)$ such that eq. (\ref{epse})
 is satisfied, i.e.

\begin{equation}
\left\{\begin{array}{l} A^+A+\epsilon C^+C={\oo}_n,\\
        \epsilon B^+B+D^+D={\oo}_m,\\
        \epsilon B^+A+D^+C={\bf 0}.\end{array}\right.\label{uu}
\end{equation}

It can be shown that the equations (\ref{uu}) also imply the equivalent
 relations
\begin{equation}
\left\{\begin{array}{l}AA^++\epsilon BB^+={\oo}_n,\\
\epsilon CC^++DD^+={\oo}_m,\\
\epsilon AC^++BD^+={\bf 0}.\end{array}\right.\label{uuu}
\end{equation}

Now we find the matrix $U$ with the property $Z'(Z_1)={\bf 0}$, i.e.
$$AZ_1+B={\bf 0}.$$

With  the first eq. (\ref{uuu}), it is obtained

$$A^+A=(\oo +\epsilon Z_1Z_1^+)^{-1}~.$$

A polar decomposition of the matrix $A$ is used
$$A=XH~,$$
where $H$ is hermitian and positive definite, while $X$ is unitary.

 The matrix $U$ has as subblocks the submatrices
\begin{equation}
\left\{\begin{array}{l} A=X(\oo +\epsilon Z_1Z_1^+)^{-1/2},\\ \label{712}
                        B=-X(\oo +\epsilon Z_1Z_1^+)^{-1/2}Z_1,\\
                        C=\epsilon X'(\oo +\epsilon Z_1^+Z_1)^{-1/2}Z^+_1,\\
                        D=X'(\oo +\epsilon Z_1^+Z_1)^{-1/2},\end{array}\right.
\end{equation}
where $X~(X')$ is a unitary $n\times n$ (resp. $m\times m$)  matrix,
irrelevant for the calculation of the distance.

Note that the representation (\ref{712}) of the matrix $U$ coincides with the
one given  by eq. (\ref{56}), with $Z$ replaced by $-Z$. This is a consequence
of the fact that the representation (\ref{56}) expresses the  transformation
 $0\rightarrow Z$. The inverse transformation  is given by equation
(\ref{zz'}) below. 

The condition ${\rm Det}\, U=1$ is verified with the Schur formulas
 (\ref{schur}) and the representation (\ref{zz}) is obtained by $Z'=Z'(Z_2)$. 

Note also that  the linear fractional transformation can also be written down
as 
\begin{eqnarray}
Z' & = & (AZ+B)(CZ+D)^{-1}=(ZB^+ -\epsilon A^+)^{-1}(C^+-\epsilon ZD^+),
\label{z'z}\\
Z  & = & (A-Z'C)^{-1}(Z'D-B)=(C^++\epsilon A^+Z')(B^+Z'+\epsilon D^+)^{-1}.
\label{zz'}
\end{eqnarray}

With these relations it is easy to show that the homographic
 transformations (\ref{omo}) leave invariant the equation (\ref{65}) of
 geodesics.
To verify the last assertion, the following relations are also needed
\begin{equation}
\begin{array}{l}\oo +\epsilon Z^+Z=(Z'^+B+\epsilon D)^{-1}
(\oo +\epsilon Z'^+Z')
(B^+Z'+\epsilon D^+)^{-1},\\
dZ=-(\epsilon Z'C-\epsilon A)^{-1}dZ'(B^+Z'+\epsilon D^+)^{-1},\\
d^2Z=(\epsilon Z'C-\epsilon A)^{-1}[2\epsilon dZ'C
(\epsilon Z'C-\epsilon A)^{-1}dZ'-d^2Z'](B^+Z'+\epsilon D^+)^{-1}.
\label{fr}
\end{array}
\end{equation}

$b)$ Now the distance between the points $Z_1=0,~Z_2=Z$ is calculated, where
 in the compact case $Z\in {\cal V}_0$.

The starting point is the formula (\ref{dist}). Using eqs. (\ref{fr}) and
 also the equation
$$\oo +\epsilon ZZ^+=(Z'C-A)^{-1}(\oo +\epsilon Z'Z'^+)
(C^+Z'^+- A^+)^{-1},$$
it is easy to verify that the infinitesimal element (\ref{dist}) is invariant
under the homographic transformations (\ref{omo}).

We fix $k=1$ in eq. (\ref{dist}). The expression (\ref{58}) of the geodesics
 with $B\rightarrow Bt$ implies
\begin{equation}
(\oo +\epsilon ZZ^+)^{-1}={\rm co}^{2}(t\sqrt {BB^+});
~(\oo +\epsilon Z^+Z)^{-1}=
{\rm co}^{2}(t\sqrt {B^+B}),
\end{equation}
and the equation (\ref{dist}) can be written down as 
\begin{equation}
ds^2={\rm Tr}\,(B^+B)dt^2=\sum |B_{ij}|^2dt^2,
\end{equation}
which shows that indeed $B_{ij}\in {\got m}$ are normal coordinates.

So, we find eq. (\ref{tr}) and the first eq. (\ref{un}):
$$d^2(0,Bt)={\rm Tr}\,(B^+B)t={\rm Tr}\,[{\rm arcta}(ZZ^+)^{1/2}]^2,$$
which implies also the other two relations in eqs. (\ref{un}) and 
(\ref{lag}).

$c)$ Finally, the representation (\ref{v}) for the matrix $V$ is obtained
after a tedious matrix calculation. We only point out the main steps.

If in eq. (\ref{zz}) it is substituted
$$Z_2-Z_1=(\oo +\epsilon Z_1Z^+_1)Z_2-Z_1(\oo +\epsilon Z^+_1Z_2),$$
then it is obtained
$$Z=(\oo +\epsilon Z_1Z^+_1)^{1/2}(\oo +\epsilon Z_2Z^+_1)^{-1}Z_2
(\oo +\epsilon Z^+_1Z_1)^{1/2}-Z_1,$$
\begin{equation}
(\oo +\epsilon Z_1Z_1^+)^{-1/2}Z=(\oo +\epsilon Z_2Z^+_1)^{-1}
Z_2(\oo +\epsilon Z^+_1Z_2)^{1/2}-(\oo +\epsilon Z_1Z^+_1)^{-1/2}Z_1.
\label{n}
\end{equation}

The representation (\ref{n}) is introduced in the auxiliary expression

$$E=(\oo +\epsilon Z_1Z_1^+)^{-1/2}ZZ^+(\oo +\epsilon Z_1Z_1^+)^{-1/2}.
\label{ee}$$
 It is obtained
$$E=F+(\oo +\epsilon Z_1Z^+_1)^{-1}Z_1Z_1^+,$$
where, finally, $F$ is brought to the form
$$F=(\oo +\epsilon Z_2Z^+_1)^{-1}(\epsilon\oo +Z_2Z^+_2)
(\oo+\epsilon Z_1Z_2^+)^{-1}
-\epsilon  \oo.$$

So, $E$ can be written down as
$$E=(\oo +\epsilon Z_2Z^+_1)^{-1}(\epsilon\oo +Z_2Z^+_2)(\oo 
+\epsilon Z_1Z^+_2)^{-1}-\epsilon (\oo +\epsilon Z_1Z^+_1)^{-1},$$
which implies
$$
 \oo +\epsilon ZZ^+=(\oo +\epsilon Z_1Z_1^+)^{1/2}(\oo +\epsilon Z_2Z^+_1)^{-1}
(\oo +\epsilon Z_2Z^+_2)(\oo +\epsilon Z_1Z^+_2)^{-1}
(\oo +\epsilon Z_1Z^+_1)^{1/2},\label{cto}
$$
i.e. the representation (\ref{v}).~\pre

Now we comment the expressions given by Proposition {\ref{pr5}}:

\begin{com}
\label{com3} The formulas giving the distance on the complex Grassmann
manifold (his noncompact dual) generalize the corresponding ones for the
Riemann sphere (resp., the disk $|z|<1$).
\end{com}
 Formula (\ref{r}) was used by Rosenfel'd,\cite{ros}
$\theta_i$ being the stationary angles here determined in Lemma 4, while the
first expression (\ref{lag}) appears  in Ch.
 II \S6 p. 69  of Ref. \onlinecite{pia} or in  Ref. \onlinecite{sieg} in the
 case of symplectic group.
 In  the last case the factor $k=4$ in formula (\ref{dist}) eliminates the
 factor $1/2$ in  the first eq. (\ref{lag}).

Now we particularize the formulas in
 Proposition 3 to the case of the Riemann sphere.

If $Z_1,~ Z_2$ belong to the same chart, then

\begin{equation}
Z=\frac {Z_1-Z_2}{1+\bar{Z}_1Z_2},~ V=\frac{|1+Z_1\bar{Z}_2|^2}{(1+|Z_1|^2)
(1+|Z_2|^2)}~,
\end{equation}
and formulas (\ref{un}) become
\begin{eqnarray}
d(Z_1,Z_2) & = & \arctan\frac{|Z_1-Z_2|}{|1+\bar{Z}_1Z_2|}
  \nonumber = 
  \arccos \frac{|1+Z_1\bar{Z}_2|}{(1+|Z_1|^2)^{1/2}(1+|Z_2|^2)^{1/2}}\\
  & = & \arcsin \frac{|Z_1-Z_2|}{(1+|Z_1|^2)^{1/2}(1+|Z_2|^2)^{1/2}}~.
\end{eqnarray}

Note that
\begin{equation}
d(Z_1,Z_2)=d_c(Z_1,Z_2)=\frac {1}{2}\theta (Z_1,Z_2)=
\arcsin \frac{\Delta (Z_1,Z_2)}{2}~,
\end{equation}
where $\theta (Z_1,Z_2)$ is the length of the arc of the great circle
 (the geodesic) joining
the points $Z_1,~Z_2\in {\gcp}^1$, while $\Delta (Z_1,Z_2)$ is the chord
length.

Eq. (\ref{lag}) in the case of Riemann sphere (respectively the disk $|Z|<1$) 
($\epsilon =1,~ \eta= i~ (\epsilon =-1,~\eta =1$)) reads
\begin{equation}
d(Z_1,Z_2)=\frac{1}{2\eta}\log\frac{1+\eta\lambda}{1-\eta\lambda}~,~ \lambda =
\frac{|Z_1-Z_2|}{|1+\epsilon Z_1\bar{Z}_2|} ~.
\end{equation}

The expression under the logarithm represents the cross-ratio
 $\{{\rm Z}_1{\rm Z}_2,{\rm MN}\}$. In the compact (noncompact) case  M and N
represents the points where the line  ${\rm Z}_1,{\rm Z}_2$ meets the absolute
(Laguerre-Cayley-Klein) (resp. the frontier $|Z|=1$).\cite{cox} 

In the case of ${\gcp}^n$, the second relation in eq. (\ref{un}) is the
elliptic hermitian Cayley distance (\ref{cayley})  expressed in non-homogeneous
 coordinates
\begin{equation}
d(Z,Z')=d_c(Z,Z')=\arccos \frac{|1+\sum Z_i\bar{Z'}_i|}{(1+\sum |Z_i|^2)^{1/2}
(1+\sum |Z'_i|^2)^{1/2}}~,
\end{equation}
while in the case of the hermitian hyperbolic space 
$SU(n,1)/S(U(n)\times U(1))$ the
corresponding  distance is the hyperbolic hermitian distance (\ref{ncay})
\begin{equation}
d(Z,Z')={\rm arccosh}\, \frac{|1-\sum Z_i\bar{Z'}_i|}{(1-\sum |Z_i|^2)^{1/2}
(1-\sum |Z'_i|^2)^{1/2}}~.~\pre
 \end{equation}

Finally, let us denote by $\delta _n$ the distance (\ref{r}) on \Gras~ and his 
noncompact dual (\ref{xn})
 and by $s_n$ the distance (\ref{55}) of the images of the points through the
 Pl\"{u}cker embedding (resp. (\ref{55n})), where both points belong to
 ${\cal V}_0$ for $X_c$. We present an elementary inequality which has a simple
geometrical meaning in the following  

\begin{com}
\label{comas}
Let $\delta_n$ and $s_n$ be defined by
\begin{eqnarray}
\label{frst1} \delta_n^2 & = & \theta_1^2+\ldots +\theta_n^2~,\\
\label{second1}\mbox{\rm co}\,s_n & = & \mbox{\rm co}\,\theta_1\cdots
 \mbox{\rm co}\,\theta_n~,
\end{eqnarray}
where  for $X_c$ $s_n, ~\theta_i\in [0,\pi /2]$. Then
\begin{equation}
\delta_n\geq s_n~.
\label{xinn1}
\end{equation}
\end{com}

{\it Proof\/}:~ We indicate the proof on $X_c=\Gras$.  Geometrically, eq. 
(\ref{xinn1})
expresses the fact that the distance on the manifold (here the Grassmannian)
is greater than the distance between the images of the points through the
 embedding (here the Pl\"{u}ckerian one). Infinitesimally
\begin{equation}
\label{infini}
d\,\delta_n=d\,s_n,
\end{equation}
as can be verified with eq. (\ref{second1})  at small stationary angles. 

We also indicate an elementary algebraic proof of (\ref{xinn1}). Firstly, eq.
 (\ref{xinn1})
is proved for $n=2$, i.e. if $\delta^2=x^2+y^2$ and $\cos s=\cos x\cos y$,
 then
$\delta\geq s$ for $x, y\in [0,\pi /2]$. Indeed, let $x=\delta\cos \theta , 
~y=\delta\sin \theta$
and let us consider the function $F(\theta)\equiv\cos s$. Then the equation
 $dF/d\theta =0$ has as unique solution in $[0,\pi /2]$ the angle
 $\theta =\pi /4$, which is a maximum. But $F(0)=F(\pi /2)=\cos\delta$ and
the inequality $\cos \delta\leq \cos s$ follows. Further the mathematical 
induction
on $n$ in eqs. (\ref{frst1})-(\ref{xinn1}) is applied.~\pre
\vspace{1cm}

\hspace{-4mm}\hspace{2mm}  ACKNOWLEDGEMENTS 
\vspace*{0.8cm}

This paper was began in 1993 and continued in 1994 at the Institute
 de Math\'ema\-tique de l'Universit\'e Paris 7 Denis Diderot. The author thanks
to Prof. Anne Boutet de Monvel and CNRS for the opportunity to work in the 
\'Equipe de Physique Math\'ematique et G\'eom\'etrie de CNRS. The author has
largely benefited of the interest and suggestions of Professor L. Boutet de
 Monvel. The author expresses his thanks to Professor
K. Teleman, Dr. B. Berceanu and Dr. C. Gheorghe for illuminating discussions.
 The author also expresses his thanks to Professors M. Berger, Th. Hangan,
 S. Kobayashi, B. Rosenfel'd and T. Sakai for correspondence and suggestions.
 The author expresses his gratitude to Prof. L. Papaloucas
for inviting him to the Department of Mathematics of the University of Athens,
Panepistemiopolis, where a part of this article was completed. The author
 thanks to the referee for the useful suggestions.

\rigla
\footnotesize


\begin{thebibliography}{199}


\bibitem{postnikov}M. Postnikov, {\it Lectures in Geometry, Semester II,
Linear Algebra and Differential Geometry} (Mir Publishers, Moscow, 
 1982), p. 105.


\bibitem{jo}C. Jordan, ``Essai sur la G\'eom\' etrie \`a n dimensions'',
 Bull. Soc. Math. France  {\bf t. lll} 103 (1875),  in 
{\it \OE uvres de C. Jordan}, Paris, Gauthier-Villars, {\bf Tome III} 79 
(1962).

\bibitem{eh}C. Ehresman,`` Sur la topologie de certain espaces homogen\`es'',
  Ann. Math. {\bf 35} 396 (1934).

\bibitem{chern}S. S. Chern, {\it Complex Manifolds without Potential
Theory} (Van Nostrand, Princeton, 1967).

\bibitem{pont}L. C. Pontrjagin, ``Charakteristiceskie tzikly differentziruemyh
mnogobrazia'',   Mat. sb. {\bf 21} 233 (1947). 

 \bibitem{bour}N. Bourbaki {\it \'El\'ements de Math\'ematique, Alg\`ebre},
Chapitres 1 \`a 3  (Hermann, Paris, 1970).

\bibitem{hodge}W. D. Hodge and D. Pedoe, {\it Methods of Algebraic Geometry},
Vol. II (Cambridge University, Cambridge, 1952).

\bibitem{hus}D. Husemoller, {\it Fibre Bundles} (Mc Graw-Hill, New York,
1966).

\bibitem{mil}J. Milnor and J. D. Stasheff, {\it Characteristic Classes},
Annals of Maths. Studies Vol. 76 (Princeton, New Jersey, 1974).

\bibitem{wong}Y. -C. Wong, ``Differential Geometry of Grassmann manifolds'',
 Proc. Nat. Acad. Sci. U.S.A. {\bf 57}  589 (1967).

\bibitem{ros}B. Rosenfel'd, ``Vnnutrenyaya geometriya mnojestva 
m-mernyh ploskastei n-mernova ellipticeskovo prostranstva'',
Izv. Akad. Nauk. SSSR, ser. mat. {\bf 5} 353 (1941).

\bibitem{rosi1}B. Rosenfel'd, {\it Mnogomernye Prostranstva} (Nauka, Moskwa,
1966).

\bibitem{rosi}B. Rosenfel'd,  {\it Neevklidovy Pronstranstva} (Nauka, Moskwa,
 1969).

\bibitem{won}Y. -C. Wong, ``Conjugate loci in Grassmann manifold'',
Bull. Am. Math. Soc.  {\bf 74} 240 (1968).

\bibitem{kn}S. Kobayashi and K. Nomizu, {\it Foundations of Differential
Geometry}, Vol II (Interscience, New York, 1969).

\bibitem{wo}Y. -C. Wong, ``A class of Schubert varieties'',  J. Diff.
Geom. {\bf 4} 37 (1970).


\bibitem{sak}T. Sakai, ``On cut loci on compact symmetric spaces'',
 Hokkaido Math. J. {\bf 6} 136 (1977).

\bibitem{kob}S. Kobayashi, ``On conjugate and cut loci'',  in
{\it Global Differential Geometry}, M.A.A. Studies in Mathematics,
Vol.  27, S. S. Chern editor (1989),  p. 140.

\bibitem{wwong}Y. -C. Wong, ``Isoclinic $n$-planes in Euclidean $2n$-space,
Clifford parallels in elliptic $(2n-1)$ space, and the Hurwitz matrix
 equations'',
 Mem. Amer. Math. Soc. {\bf 41} (1961).

\bibitem{hur}A. Hurwitz, ``\"{U}ber die Komposition der quadratischen Formen
 von
 beliebig vielen Variablen'',  Nach. v. der Ges. der. Wiss., G\"{o}tingen
(Math. Phys. K1)  309 (1898); reprinted in {\em Math. Werke}
 Bd. 2, p. 565.

\bibitem{per}A. M. Perelomov, ``Coherent states for arbitrary Lie groups'',
 Commun. Math. Phys. {\bf 26} 222 (1972). 

\bibitem{klauder} {\it Coherent states}, edited by J. R. Klauder and
 B. S. Skagerstam (Word Scientific, Singapore, 1985).

\bibitem{sb}S. Berceanu, ``The coherent states: old geometrical methods
in new quantum clothes'', preprint Bucharest, Institute of Atomic
 Physics, FT-398 (1994) and also Bielefeld preprint Nr.
 664/11/1994.

\bibitem{sb1}S. Berceanu, ``Coherent states and global differential geometry'',
paper presented at the XIII Workshop on Geometrical methods in Physics, 
Bia\l owe\.{z}a 1994,  in {\it Quantization, coherent
states and complex structures} edited by J.-P. Antoine et all,
 (Plenum, New York, 1995) p. 131 ;
``Coherent states and
 geodesics:
cut locus and conjugate locus", J. Geom. Phys. {\bf 21} 149 (1997).

\bibitem{ps}A. Pressley and G. Segal, {\it Loop Groups} (Oxford University,
 Oxford, 1986).

\bibitem{cal} E. Calabi, ``Isometric imbedding of complex manifolds'',
 Ann. of Math. {\bf 58} 1 (1953).

\bibitem{cgr}M. Cahen, S. Gutt and J. Rawnsley, ``Quantization on K\"ahler
manifolds'', ll, Trans. Math. Soc.  {\bf 337} 73 (1993).  



\bibitem{leich}K. Leichtweiss, ``Zur Riemannschen Geometrie in Grassmannschen
Mannigfaltigkeiten'',   Math. Z. {\bf 76} 331 (1961).





\bibitem{ms}M. Spera and G. Valli, ``Pl\"ucker embedding of the Hilbert space
 Grassmannian and the CAR algebra'', Russian J. Math. Phys. {\bf 2} 383 (1994).

\bibitem{gant}F. P. Gantmacher, {\it The Theory of Matrices}
(Chelsea,  New York, 1959).




\bibitem{cay}A. Cayley, ``A sixth memoir upon quantics'',
 Phil. Trans. Royal. Soc. London {\bf 149} 61 (1859),  in
{\it Collected Mathematical Papers}, (Cambridge, 1889-1898).

\bibitem{coo}J. L. Coolidge, ``Hermitian metrics'',  Ann. of  Math.
{\bf  22} 11 (1921).

\bibitem{kobi}S. Kobayashi, ``On the Geometry of bounded domains'',
 Trans. Amer. Math. Soc. {\bf 92} 267 (1959).


\bibitem{study}E. Study, ``K\"urzeste wege im complexe gebiete'',
 Math. Ann. {\bf 60} 321 (1905).

\bibitem{cool}J. P. Coolidge, {\it The Geometry of complex domain}
(Clarendon, Oxford, 1924).


\bibitem{hirz}F. Hirzebruch, {\it Topological Methods in Algebraic Geometry}
(Springer-Verlag, Berlin, 1966).



\bibitem{wu}H. H. Wu, {\it The Equidistribution Theory of Holomorphic Curves},
Annals of Maths. Studies Vol. 164 (Princeton, NJ, 1970).





\bibitem{lu}A. Lundell and S. Weingram, {\it The Topology of CW Complexes}
(Van-Nostrand-Reinhalt, Toronto, 1969).


\bibitem{sbcag}S. Berceanu and A. Gheorghe, ``On the construction of
perfect Morse functions on compact manifolds of coherent states'',
  J. Math. Phys.   {\bf 28} 2899 (1987).






\bibitem{som}D. M. Y. Sommerville, {\it An Introduction to the Geometry of n
Dimensions} (Dover,  New York, 1929).

\bibitem{hi}S. Hildebrandt, J. Jost and K. O. Widman, ``Harmonic mappings and
minimal submenifolds'',  Invent. Math. {\bf 62} 269 (1980).

\bibitem{fc}D. Fischer-Colbrie, ``Some rigidity theorems for minimal
 submenifolds of the sphere'',  Acta. Math. {\bf 145} 29 (1980).

\bibitem{wongg}Y. -C. Wong, ``Euclidean $n-$spaces in pseudo-euclidean spaces
 and
differential  Geometry of Cartan domains'',  Bull. Am. Math. Soc. {\bf 74} 409
 (1969).


\bibitem{helg}S. Helgason, {\it Differential Geometry, Lie groups and
Symmetric Spaces} (Academic, New York, 1978).

 \bibitem{sbl}S. Berceanu and L. Boutet de Monvel, ''Linear dynamical
systems, coherent state manifolds, flows and matrix Riccati equation'', 
 J. Math. Phys.  {\bf 34} 2353 (1993). 

\bibitem{wolf}J. A. Wolf, ``The action of a real semisimple Lie group on  a
complex flag manifold 1: Orbit structure and holomorphic arc components'', 
 Bull. Am. Math. Soc. {\bf 75} 1121 (1969).
\bibitem{tel}K. Teleman, ``Sur les vari\'et\'es de Grassmann'',  Bull. Math.
Soc. Sci. Math. Phys. Rep. Popul. Roum. {\bf 2(50)} 202 (1958).


\bibitem{knapp}A. W. Knapp, {\it Representation Theory of Semisimple Lie
Groups} (Princeton, NJ, 1986).




\bibitem{cr}R. Crittenden, ``Minimum and conjugate points in symmetric
 spaces'',  Canad. J. Math. {\bf 14} 320 (1962).

\bibitem{war}F. W. Warner, ``The conjugate locus of a Riemannian manifold'',
 Amer. J. Math. {\bf 87} 575 (1965).

\bibitem{wei} A. D. Weinstein, ``The cut locus and conjugate locus of a
Riemannian manifold'',  Ann. of Math. {\bf 87} 29 (1968).



\bibitem{sak1}T. Sakai, ``On the structure of cut loci in compact Riemannian
symmetric spaces'',  Math. Ann. {\bf 235} 129 (1978).

\bibitem{hua}L. K. Hua, {\it Harmonic Analysis of Functions of several Complex
Variables in the classical Domains} (Science Press, Peking, 1958). Russian
translation from Chinese (Izdat. Inostr. Lit., Moscow, 1959). English
translation from Russian Translations of Math. Monographs,
Amer. Math. Soc. Vol. 6 (Providence, NY, 1963).

\bibitem{ben}A. Ben-Israel, ``Generalised inverse matrices: a perspective of
 the work of Penrose'',  Proc. Cambridge. Philos. Soc. {\bf 100} 407 (1986).



\bibitem{ww1}J. A. Wolf, ``Geodesic spheres in Grassmann manifolds'',  Ilin. J.
Math. {\bf 7} 425 (1963); ``Elliptic spaces in Grassmann manifolds'', 
 Ilin. J. Math. {\bf 7} 447 (1963).


\bibitem{pia} I. I. Pjatetskii-Shapiro, {\it Geometry of classical Domains and
Theory of Automorphic Functions}, (Gasudarstvennoe izdatelistvo
fiziko-matematicheskoi literatury, Moskwa, 1961); French translation from
Russian (Dunod, Paris, 1966).

\bibitem{sieg}C. L. Siegel, ``Symplectic Geometry'',  Amer. J. Math. {\bf 65} 
1 (1943).

\bibitem{cox}H. S. M. Coxeter, {\it Non-Euclidean Geometry}
(Toronto, Toronto, 1942).



\end{thebibliography}
\end{document}